\documentclass[journal]{IEEEtran}
\ifCLASSINFOpdf
\else
\fi
\usepackage{epsfig}\usepackage{epsfig}
\usepackage{amssymb}
\usepackage{amsthm}
\usepackage{amsmath}
\usepackage{cite}
\usepackage{diagbox}
\usepackage{multirow}
\usepackage{mathrsfs}
\usepackage{cases}
\usepackage{algorithm}
\usepackage{multirow}
\usepackage{makecell}
\usepackage{algorithmic}
\usepackage{subeqnarray}
\usepackage{subfigure}
\usepackage{pifont,bbm,amsfonts,mathrsfs,amsmath,stmaryrd,amssymb}
\usepackage{amscd,graphicx,array,dsfont,texdraw,tikz,footnote,verbatim}%
\usepackage[normalem]{ulem}
\usepackage{mathtools}
\usepackage{bbm,bbding,amssymb,pifont,mathrsfs,amsfonts,graphicx,mathtools,color}%
\usepackage{amscd,array,enumerate,dsfont,texdraw,tikz,multicol,bbm}

\usetikzlibrary{arrows,shapes,snakes,shadows,positioning,automata,patterns}
\usetikzlibrary{trees,decorations.pathmorphing,decorations.markings}

\usepackage{tikz}
\usepackage{pgf}
\usepackage{amsmath, amsfonts, amssymb, amsthm}
\usepackage{pgfplots}
\pgfplotsset{compat=newest}
\usepackage{xxcolor}
\usepackage{longtable}
\usepackage{array}

\usepackage{booktabs}
\usetikzlibrary{arrows, decorations.markings}
\usetikzlibrary{arrows,shadows,petri}
\newcommand{\beq}{\begin{equation}}
\newcommand{\eeq}{\end{equation}}
\newcommand{\bqa}{\begin{eqnarray}}
\newcommand{\eqa}{\end{eqnarray}}

\definecolor{maroon}{rgb}{0.7,0,0}

\definecolor{ngreen}{rgb}{0.3,0.7,0.3}

\definecolor{golden}{rgb}{0.8,0.6,0.1}

\DeclareMathOperator*{\argmin}{arg\,min}

\newtheorem{theorem}{\indent Theorem}
\newtheorem{lemma}{\indent Lemma}

\newtheorem{definition}{\indent Definition}

\newtheorem{myremark}{\indent Remark}

\newenvironment{remark}{\begin{myremark}\normalfont}
	{\end{myremark}}

\begin{document}
\title{Defense for Advanced Persistent Threat with Inadvertent or Malicious Insider Threats}

\author{Ziqin Chen, Guanpu Chen and Yiguang Hong, \textit{Fellow, IEEE}
\thanks{This work was supported by the National Natural Science Foundation of China (No. 62173250), by Shanghai Municipal Science and Technology Major Project (No. 2021SHZDZX0100), and by the China Postdoctoral Science Foundation (No. 2021M702481).}
%by National Key R\&D Program of China under Grants 2018YFE0105000 and 2018YFB1305304, by Shanghai Municipal Science and Technology Major Project under Grant 2021SHZDZX0100, by Shanghai Municipal Science and Technology Fundamental Project under Grant 21JC1405400, and by the Shanghai Municipal Commission of Science and Technology under Grants 1951113210 and 19511132101.}
\thanks{Z. Chen is with the Department of Control Science and Engineering \& Shanghai Research Institute for Intelligent Autonomous Systems, Tongji University, Shanghai, 201804, China (e-mail: cxq0915@tongji.edu.cn).} 
\thanks{G. Chen is with JD Explore Academy, Beijing, 100176, China  (e-mail:  chengp@amss.ac.cn).}
\thanks{Y. Hong is with Department of Control Science and Engineering \& Shanghai Research Institute for Intelligent Autonomous Systems, Tongji University,
Shanghai, 201804, China, and is also with Key Laboratory of Systems and Control, Academy of Mathematics and Systems Science, Beijing, 100190,
China (e-mail: yghong@iss.ac.cn).}}%<-this stops a space
\maketitle

%{}

% . Modern society is highly vulnerable to advanced persistent threats (APTs).When APTs are coupled with insider threats, the situation becomes complicated. According to a different intrinsic motivation of an insider, it should be distinguished between an inadvertent and a malicious insider.

\begin{abstract}
In this paper, we propose a game-theoretical framework to investigate advanced persistent threat problems with two types of insider threats: malicious and inadvertent. Within this framework, a unified three-player game is established and Nash equilibria are obtained in response to different insiders. By analyzing Nash equilibria, we provide quantitative solutions to the advanced persistent threat problems with insider threats. Furthermore, optimal defense strategy and defender's cost comparisons between two insider threats have been performed. The findings suggest that the defender should employ more active defense strategies against inadvertent insider threats than against malicious insider threats, despite the fact that malicious insider threats cost the defender more. Our theoretical analysis is validated by numerical results, including an additional examination of the conditions of the risky strategies adopted by different insiders. This may help the defender in determining monitoring intensities and defensive strategies. 
\end{abstract}


\IEEEpeerreviewmaketitle

% Unlike traditional cyberattackers, APTs are executed by well-resourced and well-organized attackers and characterize the continuous interplays of advanced attacks on the system resources over a long-time.

\section{Introduction}
%APT简单介绍.
\IEEEPARstart{A}{d}vanced Persistent Threat (APT) as a new type of cyber attack has posed a severe threat to nation-states and enterprise networks~\cite{applying1,applying2}. The APT is executed by a well-resourced and well-organized entity in order to steal sensitive data
covertly and on a long-term basis from the target organization. Since ``Stuxnet" attacked Iran's nuclear power plants~\cite{example1}, APTs, such as Ocean Lotus, Duqu, and Flame have attracted much attention~\cite{example2}. Currently, one practical approach for mitigating the majority of APTs is the game-theoretical method~\cite{FlipIt,securegame2,securegame3,securegame4,securegame5,securegame6}, which provides insights on feasible defense decision-making through mutual strategic behavior analysis. Traditional APT problems centered on defending against external attacks with ignoring the presence of insiders.

%引出内奸威胁以及处理的难点,处理的办法，但是是两体博弈，三体更为复杂。
However, statistics showed that the damage
caused by insiders is even more serious than outsider attackers~\cite{insiderthreat}, because the insider has privileged access to sensitive resources and system accounts. Hence, both academia and industry have paid increasing attention to insider threats in recent years~\cite{insiderchallenge,insiderpeople,insiderkuangjia}. To deal with insider threats, several approaches have been proposed, including system dynamic approaches~\cite{insiderdynamicapproachs}, machine learning based approaches~\cite{insidermachinelearningapproaches} and game-theoretic approaches~\cite{insidergameapproaches1,insidergameapproaches2,insidergameapproaches3}. Notably, the approaches presented above were primarily focused on a two-player game involving the defender and the insider.

It is worth mentioning that game-theoretic approaches to quantitatively analyzing the APT problem with insider threats (APT-I problem) are limited. Some of the reasons are as follows. On the one hand, the game model is more difficult to establish due to the interconnection among the defender, the attacker, and the insider than in two-player games of the defender and the insider. On the other hand, the coexistence of competition and cooperation makes the traditional zero-sum game method given in~\cite{insidergameapproaches1} inapplicable. To address this issue, some efforts have been made in~\cite{Feng,Hu,SCADA, dianwang}. For example, \cite{Feng} modeled the APT-I problem as a three-player timing game, in which the defender chooses appropriate time points for investing defense resources in order to minimize the APT-I impact at the lowest cost. Moreover, ~\cite{Hu} further captured the APT's conspicuous features that the long period of completed attacks with a two-layer game. Specifically, based on a differential equation for the evolution of the fraction of the compromised organization, the defend-attack between the defender and the outside attacker, as well as information trading among insiders were modeled, respectively. However, the role and impact of the insider on the intrusion process of the attack on the system were not considered in most of existing results.

Actually, the three-player games mentioned in~\cite{Feng,Hu,SCADA, dianwang} only considered the so called ``malicious insider"~\cite{insiderkuangjia}, who is motivated to collaborate with outside attackers for personal, financial, or revenge reasons. In fact, another type of insider is known as ``inadvertent insiders"~\cite{inadvertent1,inadvertent2,inadvertent3}, who are individuals without malicious intent, but their actions or inactions can jeopardize the organization's assets and operations, causing harm or increasing the likelihood of future harm to the organization's confidentiality, integrity, and availability. For example, employees accidentally leak sensitive data on social media, lose work devices, and fall victim to phishing and other disguised malware attacks~\cite{insiderkuangjia,inadvertent1,inadvertent2}. To our knowledge, no game-theoretic approaches to the APT-I problem have taken into account inadvertent insider threats. The inadvertent insider's strategies are hard to predict due to their inadvertent nature. However, the strategies of the inadvertent insider are strongly coupled with those of the defender and attacker, determining the distribution of Nash equilibria. As a result, analyzing the Nash equilibria of a three-player APT game with an inadvertent insider becomes more complicated. 

The objective of this paper is to distinguish between two types of insider threats in the APT-I problem and establish a game framework to answer the following three questions: i) What is the optimal defense strategy for dealing with the APT-I problem with inadvertent insider threats or malicious threats? ii) What are the differences in optimal defense strategies and the defender's costs between the two previously mentioned insider threats? In addition, iii) What strategies should be used to achieve a feasible defense if the defender is uncertain about the types of insider threats? The contributions are summarized as follows.
\begin{itemize}
	\item This is the first work to establish a game-theoretic framework for the APT problem, involved with two types of insiders, malicious and inadvertent. Within this framework, as the risk budget mechanism is introduced to provide an incentive for the inadvertent insider, the APT-I problem can be uniformly formulated as a three-player game played among a defender, an attacker, and an insider. Our formulation includes the two-player games in~\cite{insidergameapproaches1,insidergameapproaches2,insidergameapproaches3} and the model with only malicious insiders~\cite{Feng,Hu,SCADA, dianwang}.

	\item To deal with the two types of insider threats, we obtain the feasible Nash equilibria of the corresponding APT-I game in order to provide a quantitative analysis under insider threats. It is worth noting that the optimal defense strategies for the defender when dealing with malicious insider threats are not always identical to those for inadvertent insider threats, due to an inherent difference between the two. Furthermore, the malicious insider causes more cost to the defender than the inadvertent, while the optimal defense strategies under the malicious threat are more worse than those under the inadvertent threat.
	
	\item  The Nash equilibria and the corresponding defender's costs are validated via numerical investigations. Furthermore, we examine the conditions of risky strategies taken by two types of insiders and find that an inadvertent insider is more likely to lead to a risky strategy. Hence, if insider threats are unknown, the defender can employ a relatively high level of monitoring and defense.
\end{itemize}

The remainder of this paper is organized as follows. An APT-I problem and a three-player APT-I game are presented in Section~\ref{problemstatement}. Detailed Nash equilibrium analysis is presented in Section~\ref{neseeking}. The discussion on optimal defense strategies and the defender's costs are summarized in Section~\ref{Discussion}. The result verification via extensive numerical is investigated in Section~\ref{Example}. Finally, concluding remarks are given in Section~\ref{Conclusion}.

\section{Problem Formulation}\label{problemstatement}
\begin{figure}
	\centering
	\includegraphics[height=5cm,width=9cm]{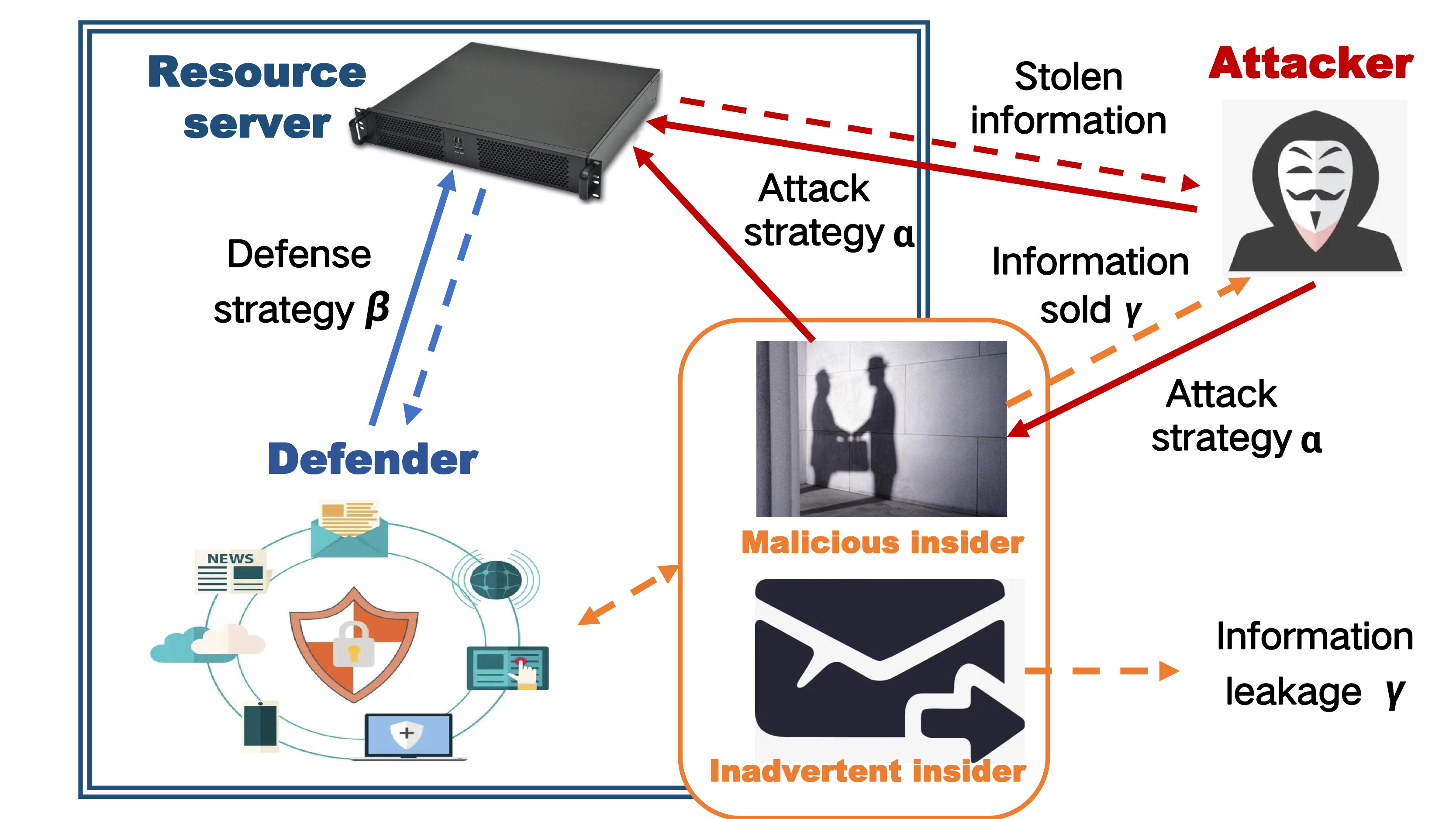}
	\caption{The diagram of an APT-I model. A single-arrowed dashed line represents the flow of privacy data from organization's resources to participants. A single-arrowed solid line represents each participant's direct action. The blue double line frame represents detection systems.}\label{fignetwork}
\end{figure}
Consider an organization with facing joint threats from an attacker and an insider, as illustrated in Fig. 1. This organization is made up of a resource server, a defender, and an insider. The resource server is used to store sensitive data, and the defender and the insider are authorized to access it for business requirements. Insiders are classified as malicious or inadvertent. A malicious insider may sell information in order to provide an attacker with back-door access. An inadvertent insider has no incentive to cooperate with an outside attacker, whose unintentional risky behavior may result in information leakage. Once the attacker has obtained information from the malicious insider, he can evade detection systems (such as firewalls, authentication controls, and intrusion prevention systems) and rapidly invade the organization to perform malicious actions, such as monitoring sensitive operations, stealing private information, and injecting ransomware. To summarize, this paper is committed to solving the following problem.

\textbf{APT-I Problem:} Consider an APT-affected organization with insider threats. How should the defender's defense strategies be allocated to minimize the system's continuing losses in response to inadvertent or malicious insider threats?

\begin{remark}
Defense strategies include user double identifications, network security monitors and employee training~\cite{applying2,inadvertent3}. In this sense, active defense strategies can be thought of as increasing the frequency of the preceding strategies. In fact, defense strategies in most organizations, such as companies, are constrained by limited computational resources. Hence, the APT-I problem is concerned with allocating finite defense strategies for mitigating the APT with insider threats.
\end{remark}

To solve the APT-I problem, we develop a mathematical model that characterizes compromised resources. In the subsequent, we refer to each participant as a \textit{player} and every feasible player's action as a \textit{strategy}.

\subsection{The evolution model of resource states}\label{systemmodel}
The organization's resources evolve dynamically as a result of the long-term persistence of the APT with insider threats. Total system resources are normalized to the value of $1$. Define $x(t)\in[0,1]$ as the fraction of compromised resources at time $t$. Correspondingly, $1-x(t)$ represents the percentage of resources under the defender's protection. The fraction of compromised resources is evolving by
\begin{equation}
	\dot{x}(t)=\alpha(1-x(t))-\beta x(t),\label{state}
\end{equation}
where attack strategies are denoted by $\alpha\in(0,1]$ and defense strategies are denoted by $\beta\in(0,1]$. At time $t$, the term $\alpha(1-x(t))$ represents the percentage of resources sized by the attacker. The term $\beta x(t)$ refers to the percentage of resources recaptured by the defender. 

Because the resource state $x(t)$ cannot be transmitted in real time to players, non-adaptive strategies are considered. In other words, strategies $\alpha$, $\beta$ and $\gamma$ do not change over time and were pre-configured prior to the organization's deployment. Non-adaptive strategies can be thought of as programmed attacks, the defender's routine security examinations, and the insider's fixed amount of information leakage. They have been proved to be a pair of optimal response strategies in the periodic APT response problem~\cite{Feng}.

Based on the evolution model \eqref{state}, the resource state at time $t$ satisfies
\begin{equation}
	x(t)=\frac{\alpha}{\alpha+\beta}\left(1-e^{-(\alpha+\beta) t}\right),\label{statesolution}
\end{equation}
which implies that organization's resources are fully under the defender's control at the beginning, i.e., $x(0)=0$. As time goes to infinity, the resource state converges to $\frac{\alpha}{\alpha+\beta}$ at a rate of $e^{-(\alpha+\beta)t}$, indicating two facts: 1) Under a sustained attack launched by an attacker, the fraction of compromised resources grows. However, due to the defense strategies adopted by a defender, the organization's resources cannot be fully invaded. 2) A large $\beta$ results in a small limit point $\frac{\alpha}{\alpha+\beta}$, implying that a active defense strategy reduces the fraction of compromised resources significantly.

%Although the insider's strategy $\gamma$ does not determine the resource state directly, it affects defense and attack strategies, which are related to resource states. 

\subsection{The APT-I game model}\label{APTImodel}
This subsection establishes a unified APT-I game framework. We start by characterizing each player's objective function in the APT-I problem.

\textbf{Attacker:} The outside attacker is dedicated to invading and occupying as many resources as possible over a long period of time at the lowest cost of attack. The attacker's cost function has three components: i) the risk of being detected by the defender, which is reduced by a factor of $(1-\gamma)^2$ with the amount of purchased information $\gamma\in[0,1]$ from the malicious insider; ii) the cost of launching attacks, which is related to attack strategies $\alpha\in(0,1]$; and iii) the cost of paying for the malicious insider. We describe the attacker's average cost function in the time horizon $[0,\infty)$ as follows,
\begin{flalign}
	J_{A}(\alpha,\beta,\gamma)&=\lim_{T\rightarrow \infty}\frac{1}{T}\int_{0}^{T}p_{A}(1-\gamma)^2(1-x(t))^2\nonumber\\
	&\quad 
	\quad+q_{A}\alpha^2(1-x(t))^2+\gamma^2(1-x(t))^2dt,\label{costattack}
\end{flalign}
where $p_{A}\in(0,1)$ is a unit risk coefficient of being detected and $q_{A}>0$ is a weight representing each attack cost. The attacker's average cost function~\eqref{costattack} develops that in the two-layer game~\cite{Hu}, in which the insider's effects $p_{A}(1-\gamma)^2(1-x(t))^2$ and $\gamma^2(1-x(t))^2$ are not taken into account.

\textbf{Defender:} The defender tries his best to minimize the damage caused by the outside attacker and the insider, whose cost function is composed of three components: i) the damage from compromised resources; ii) the operational cost of launching a defense, which is related to defense strategies $\beta\in(0,1]$; and iii) the damage from insider threats. In the time horizon $[0,\infty)$, the defender's average cost function is described as follows,
\begin{flalign}
	J_{D}(\alpha,\beta,\gamma)&=	\lim_{T\rightarrow\infty}\frac{1}{T}\int_{0}^{T}p_{D}(1-\gamma)^2x^2(t)\nonumber\\
	&\quad \quad+q_{D}\beta^2x^2(t)+\gamma^2x^2(t)dt,\label{costdefend}
\end{flalign} 
where $p_{D}\in(0,1)$ is the unit cost coefficient for compromised resources and $q_{D}>0$ is the weight standing for each defense cost. Note that the damage from compromised resources $p_{D}(1-\gamma)^2x^{2}(t)$ is reduced by a factor of $(1-\gamma)^2$ because the cost caused by an inadvertent insider $\gamma^2x^2(t)$ has been considered in the last term of the right hind side of \eqref{costdefend}. 

\begin{remark}
According to~\eqref{costdefend}, there is a positive correlation between the fraction of compromised resources $x(t)$ and the defender's average cost function $J_{D}$. Hence, minimizing the system's average loss is equivalent to finding best response strategies $\beta^*$ to minimize $J_{D}$.
\end{remark}

\textbf{Malicious Insider:} The malicious insider seeks to maximize his profits, whose profit function consists of three parts: i) the benefit from secure resources; ii) the monetary gain from trading information; and iii) the risk of being caught by a defender. To be specific, the profit function in the time horizon $[0,\infty)$ is described as follows,
\begin{flalign}
	J_{I}(\alpha,\beta,\gamma)&=\lim_{T\rightarrow \infty}\frac{1}{T}\int_{0}^{T}p_{I}(1-x(t))^2+\gamma^{2}(1-x(t))^2\nonumber\\
	&\quad \quad -(q_{I}\gamma+\frac{1}{2}\gamma^2)dt,\label{costinside}
\end{flalign} 
where $p_{I}\in(0,1)$ represents the proportion of insiders in the organization and $q_{I}\geq0$ represents the risk coefficient. \eqref{costinside} includes the risk term $q_{I}\gamma+\frac{1}{2}\gamma^2$ in addition to the insider's profit function described in the three-player game~\cite{Feng}. This risk term configuration implies that the malicious insider can not abuse his privileges without incurring any cost.

\textbf{Inadvertent Insider:} The cost of inadvertent insider's risky strategies is borne by the organization rather than the inadvertent insider, so there is no incentive for the inadvertent insider to avoid risky behaviors and worry about any potential risk caused by his actions and strategies~\cite{inadvertent1,inadvertent2,inadvertent3}. To encourage the inadvertent insider to self-manage his risks while discouraging him from his risky actions, we propose providing all employees, including the inadvertent insider, with a risky budget mechanism. This mechanism takes the same form as~\eqref{costinside} and includes three components: i) If an employee behaves prudently, the organization can reward him, as shown by $p_{I}(1-x(t))^2$; ii) If an employee engages in a risky strategy, he may face penalties, which is depicted by $-(q_{I}\gamma+\frac{1}{2}\gamma^2)$; iii) Of course, his risky strategy may bring some benefits (obtaining entertainments or relieving pressures), as shown by $\gamma^{2}(1-x(t))^2$. Under the risk budget mechanism, all employees evaluate their risk before choosing risky strategies, such as sending emails that violate data protection regulations or clicking on a malicious link in a received email~\cite{inadvertent3}. Hence, this mechanism offers incentives to an inadvertent insider while shifting the cost of risk away from the organization.

\begin{remark}
Although the risk budget mechanism's form is intended to be the same as the malicious insider's profit function~\eqref{costinside}, the optimal strategies of all players are not the same because there is no cooperation between the inadvertent insider and the outside attacker. 
\end{remark}

The balance among all players' optimal strategies is equivalent to the following Nash equilibrium solution concept based on the above objective functions~\eqref{costattack}-\eqref{costinside}.
\begin{definition}\label{defne}(Nash Equilibrium)
A strategy pair $(\alpha^*,\beta^*,\gamma^*)$ is defined as a Nash equilibrium if 
	\begin{eqnarray}
		&&J_{A}\left(\alpha^{*},\beta^{*},\gamma^{*}\right) \leq J_{A}\left(\alpha,\beta^{*},\gamma^{*}\right), \nonumber\\
		&&J_{D}\left(\alpha^{*},\beta^{*},\gamma^{*}\right) \leq J_{D}\left(\alpha^{*},\beta,\gamma^{*}\right), \nonumber\\
		&&J_{I}\left(\alpha^{*},\beta^{*},\gamma^{*}\right) \geq J_{I}\left(\alpha^{*},\beta^{*},\gamma\right).\nonumber
	\end{eqnarray}
\end{definition}
It is known that the Nash equilibrium provides a credible prediction of the attacker’s and insider's moves, identifying the optimal countermeasures to the attacker’s and insider's strategies. 

The APT-I problem can be uniformly expressed as the unified three-player game shown below.

\textbf{APT-I Game: } Consider three players in the APT-I problem: the attacker with attack strategies $\alpha$, the defender with defense strategies $\beta$, and the insider (inadvertent or the malicious insider) with strategies $\gamma$. The attacker's and the defender's goals are to minimize their costs $J_{A}$ and $J_{D}$. The insider is devoted to maximizing his profit $J_{I}$.

It should be noted that, in the worst-case scenario, the attacker may have known the exact form of $J_{A}(\alpha,\beta,\gamma)$ through meticulous reconnaissances, and may try to minimize it. In this sense, the Nash equilibrium strategies $\beta^*$ are the organization's optimal defense strategies, at least in the worst-case scenario. As a result, if $(\alpha^*,\beta^*,\gamma^*)$ is a Nash equilibrium of the ATP-I game, $\beta^*$ is an acceptable solution to the ATP-I problem. 

\section{Nash Equilibrium Analysis}\label{neseeking}
To thoroughly solve the APT-I problem and compare optimal defense strategies under two types of insider threats, we compute Nash equilibria under the following four different scenarios:
\begin{itemize}
\item Insider threats are known to the defender.

(A) There is a malicious insider in the organization;

(B) There is an inadvertent insider in the organization.

\item Insider threats are unknown to the defender.

(C) There is a malicious insider in the organization;
	
(D)  There is an inadvertent insider in the organization.
\end{itemize}
Scenarios A and B are intended to investigate the optimal defense strategies for dealing with known insider threats. Scenarios C and D provide additional optimal defense strategies for unknown threats.

\subsection{Insider threats are known}\label{defenderyes}
In this case, the defender is aware of the types of insiders. Moreover, the defender has access to insiders' strategies because the insider acts based on his privileges and accesses, and the defender is aware of his accesses. 

\subsubsection{There is a malicious insider}
In Scenario A, the APT-I game depicts the conflict between a defender and both an attacker and a malicious insider, as well as the cooperation between the attacker and the malicious insider. In terms of each player's objective functions \eqref{costattack}-\eqref{costinside}, the best response strategies of each player in Scenario A are calculated in the next lemma, whose proof is in Appendix~\ref{proofbest}.
\begin{lemma}
Considering the APT-I game with a known malicious insider, the best response strategies of the defender, attacker and malicious insider are 
\begin{flalign}
	&\alpha^{*}= \begin{cases}\frac{p_{A}(1-\gamma)^2+\gamma^2}{q_{A} \cdot \beta}, & \frac{p_{A}(1-\gamma)^2+\gamma^2}{q_{A} \cdot \beta} < 1, \\ 1, & \frac{p_{A}(1-\gamma)^2+\gamma^2}{q_{A} \cdot \beta} \geq 1.\end{cases}\label{A1}\\
	&\beta^{*}= \begin{cases}\frac{p_{D}(1-\gamma)^2+\gamma^2}{q_{D} \cdot \alpha}, & \frac{p_{D}(1-\gamma)^2+\gamma^2}{q_{D} \cdot \alpha}< 1,\\ 1, & \frac{p_{D}(1-\gamma)^2+\gamma^2}{q_{D} \cdot \alpha} \geq1.\end{cases}\label{A2}  \\
	&\gamma^{*}= \begin{cases}1, & (\frac{\beta}{\alpha+\beta})^2-q_{I}>\frac{1}{2}, \\ 1~\text{or}~~0, & (\frac{\beta}{\alpha+\beta})^2-q_{I}=\frac{1}{2}, \\ 0, & (\frac{\beta}{\alpha+\beta})^2-q_{I}<\frac{1}{2}. \end{cases}\label{A3}
\end{flalign}
\end{lemma}
According to Lemma 1, all players' strategies influence each other's best response strategy in Scenario A. Particularly, $\argmin_{\gamma\in[0,1]}J_{I}(\alpha,\beta,\gamma)\subseteq \{0,1\}$ leads to $\gamma^*\in\{0,1\}$, which means that the malicious insider can only sell all of the APT detection information or none of it at all.

\subsubsection{There is an inadvertent insider} 
In Scenario B, there is no cooperation between the attacker and the inadvertent insider because the inadvertent insider has no intention of providing the attacker with back-door access or useful information. Hence, the attacker's average cost function over the time horizon $[0,\infty)$ is described as $J_{A}(\alpha,\beta,0)$. By the similar method used in Lemma 1, the best response strategies of the defender, attacker, and inadvertent insider in terms of $J_{A}(\alpha,\beta,0)$,~\eqref{costdefend}, and~\eqref{costinside} can be obtained by
\begin{flalign}
	&\alpha^{*}= \begin{cases}\frac{p_{A}}{q_{A} \cdot \beta}, & \frac{p_{A}}{q_{A} \cdot \beta} < 1, \\ 1, & \frac{p_{A}}{q_{A} \cdot \beta} \geq 1.\end{cases}\label{B1}\\
	&\beta^{*}= \begin{cases}\frac{p_{D}(1-\gamma)^2+\gamma^2}{q_{D} \cdot \alpha}, & \frac{p_{D}(1-\gamma)^2+\gamma^2}{q_{D} \cdot \alpha}< 1,\\ 1, & \frac{p_{D}(1-\gamma)^2+\gamma^2}{q_{D} \cdot \alpha} \geq1.\end{cases}\label{B2}  \\
	&\gamma^{*}= \begin{cases}1, & (\frac{\beta}{\alpha+\beta})^2-q_{I}>\frac{1}{2}, \\ 1~\text{or}~~0, & (\frac{\beta}{\alpha+\beta})^2-q_{I}=\frac{1}{2}, \\ 0, & (\frac{\beta}{\alpha+\beta})^2-q_{I}<\frac{1}{2}. \end{cases}\label{B3}
\end{flalign}

According to~\eqref{B1}-\eqref{B3}, the inadvertent insider's strategies influence the defender's best response strategies, which further influence the attacker's best response strategies.

In addition, the inadvertent insider prefers the risk-averse strategy $\gamma^*=0$ with a relatively large penalty coefficient $q_{I}>(\frac{\beta}{\alpha+\beta})^2-\frac{1}{2}$ on the basis of \eqref{B3}. It guides the defender to suppress risky behaviors of inadvertent insiders by setting a large $q_{I}$. 

Based on~\eqref{A1}-\eqref{B3}, we provide feasible Nash equilibria of the APT-I game with two types of known insiders in Theorem 1, whose proof is in Appendix~\ref{Theorem1}.

\begin{theorem}
The Nash equilibrium strategies $\beta^*_{j},~j\in\{A,B\}$ related to $\gamma^*=1$ in Table~\ref{table1} and those related to $\gamma^*=0$ in Table~\ref{table3} are acceptable solutions to the APT-I problem with known insider threats.
\end{theorem}

\begin{table*}
	\renewcommand\arraystretch{1.8}
	\caption{\label{table1} Different Nash equilibria in Scenarios A and B}
	\begin{tabular}{|c|c|c|c|c|c|c|}
		\Xhline{1pt} 
		\multicolumn{7}{|c|}{ \textbf{Insider threats are known} } \\
		\cline{1-7} 
		\multirow{2}{*}{\textbf{Kinds}} &\multicolumn{3}{c|}{ \textbf{Malicious insider} } &
		\multicolumn{3}{c|}{ \textbf{Inadvertent insider} }   \\
		\cline{2-7} 
		& \multicolumn{2}{c|}{ \textbf{System configurations} } &
		\multirow{1}{*}{\textbf{Nash equilibria}}  & \multicolumn{2}{c|}{ \textbf{System configurations} } &
		\multirow{1}{*}{\textbf{Nash equilibria}} \\
		\cline{1-7} 
		\multirow{2}{*}{\textbf{First}} & \multirow{2}{*}{$\frac{r_{A}}{p_{A}}=\frac{r_{D}}{p_{D}}\leq1$} & \multirow{2}{*}{$q_{I}\leq (\frac{1}{\frac{r_{A}}{p_{A}}+1})^2-\frac{1}{2}$} &  \multirow{1}{*}{$\left(\frac{r_{D}}{p_{D}\beta^*_{A}},\beta^*_{A},1\right),$} 	& \multirow{2}{*}{$r_{A}=\frac{r_{D}}{p_{D}}\leq 1$} & \multirow{2}{*}{$q_{I}\leq (\frac{1}{r_{A}+1})^2-\frac{1}{2}$} & \multirow{1}{*}{$\left(\frac{r_{D}}{p_{D}\beta^*_{B}},\beta^*_{B},1\right)$,}\\
		&&&  with $\beta^*_{A}\in[\varepsilon\sqrt{p_{D}},1]$ &&&with $\beta^*_{B}\in[\varepsilon\sqrt{p_{D}},1]$ \\
		\cline{1-7} 
		\textbf{Second}	& $\frac{r_{A}}{p_{A}}\leq1,~\frac{r_{A}}{p_{A}}<\frac{r_{D}}{p_{D}}$ & $q_{I}\leq(\frac{1}{\frac{r_{A}}{p_{A}}+1})^2-\frac{1}{2}$ &  $\left(\frac{r_{A}}{p_{A}},1,1\right)$ & $r_{A}\leq1,~r_{A}<\frac{r_{D}}{p_{D}}$ & $q_{I}\leq (\frac{1}{r_{A}+1})^2-\frac{1}{2}$ & $\left(r_{A},1,1\right)$   \\
		\Xhline{1pt} 
	\end{tabular}
\end{table*}

As shown in Table~\ref{table1}, there are two kinds of Nash equilibria under various system configurations. The first kind belongs to a continuous set, whereas the second kind is a single point. Thus, different system configurations correspond to various Nash equilibria. On the other hand, the optimal defense strategies $\beta^*$ for the defender when dealing with malicious insider threats are not always identical to those for inadvertent insider threats, which implies an inherent difference between the two. From this perspective, using different defense strategies for different types of insider threats can improve the targeting and effectiveness of the APT-I defense while also conserving limited defense strategies.

\subsection{Insider threats are unknown}\label{defenderno}
In this case, the defender finds no insiders, so his average cost function over the time horizon $[0,\infty)$ is denoted by $J_{D}(\alpha,\beta,0)$. It should be noted that even if the defender knows nothing about insiders, the optimal defense strategies for different types of insider threats are not always the same. The reason is that different types of insider threats have different effects on the attacker's optimal strategies, which in turn influence the defender's optimal defense strategies. 

\subsubsection{There is a malicious insider} In Scenario C, $J_{D}(\alpha,\beta,0)$, \eqref{costattack}, and \eqref{costinside} determine the best response strategies of three players, which are calculated by using the similar method of Lemma 1 to obtain
\begin{flalign}
	&\alpha^{*}= \begin{cases}\frac{p_{A}(1-\gamma)^2+\gamma^2}{q_{A} \cdot \beta}, & \frac{p_{A}(1-\gamma)^2+\gamma^2}{q_{A} \cdot \beta} < 1, \\ 1, & \frac{p_{A}(1-\gamma)^2+\gamma^2}{q_{A} \cdot \beta} \geq 1.\end{cases}\label{C1}\\
	&\beta^{*}= \begin{cases}\frac{p_{D}}{q_{D} \cdot \alpha}, & \frac{p_{D}}{q_{D}\alpha}< 1,\\ 1, & \frac{p_{D}}{q_{D}\alpha} \geq1.\end{cases}\label{C2}  \\
	&\gamma^{*}= \begin{cases}1, & (\frac{\beta}{\alpha+\beta})^2-q_{I}>\frac{1}{2}, \\ 1~\text{or}~~0, & (\frac{\beta}{\alpha+\beta})^2-q_{I}=\frac{1}{2}, \\ 0, & (\frac{\beta}{\alpha+\beta})^2-q_{I}<\frac{1}{2}. \end{cases}\label{C3}
\end{flalign}

According to~\eqref{C1}-\eqref{C3}, the insider's strategies influence the attacker's best response strategies, which further influence the defender's best response strategies.

\subsubsection{There is an inadvertent insider}
In Scenario D, $J_{D}(\alpha,\beta,0)$, $J_{A}(\alpha,\beta,0)$, and \eqref{costinside} determine the best response strategies of three players, which are computed in the same way as Lemma 1 to obtain
\begin{eqnarray}
	&&\alpha^{*}= \begin{cases}\frac{p_{A}}{q_{A} \cdot \beta} & \frac{p_{A}}{q_{A} \cdot \beta}\leq 1, \\ 1 & \frac{p_{A}}{q_{A} \cdot \beta} > 1. \end{cases}\label{D1}\\
	&&\beta^{*}= \begin{cases}\frac{p_{D}}{q_{D} \cdot \alpha} & \frac{p_{D}}{q_{D} \cdot \alpha}\leq 1, \\ 1 & \frac{p_{D}}{q_{D} \cdot \alpha}> 1.\end{cases}\label{D2}\\
	&&\gamma^{*}= \begin{cases}1, & (\frac{\beta}{\alpha+\beta})^2-q_{I}>\frac{1}{2}, \\ 1~\text{or}~~0, & (\frac{\beta}{\alpha+\beta})^2-q_{I}=\frac{1}{2}, \\ 0, & (\frac{\beta}{\alpha+\beta})^2-q_{I}<\frac{1}{2}. \end{cases}\label{D3}
\end{eqnarray}

Observing from~\eqref{D1} and~\eqref{D2}, the strategic interaction between the defender and the attacker is identical to that in the APT game. Hence, our formulation includes the results of the two-player APT game in~\cite{Hu}.

Based on~\eqref{C1}-\eqref{D3}, feasible Nash equilibria of the APT-I game with two types of unknown insiders are provided in Theorem 2, whose proof is in Appendix~\ref{Theorem1}. For notation simplicity, some parameters are denoted as $r_{A}=\frac{p_{A}}{q_{A}}$, $r_{D}=\frac{p_{D}}{q_{D}}$, and
$$\varepsilon=\frac{1}{\sqrt{\frac{1}{r_{D}}\left(\frac{1}{\sqrt{1/2+q_{I}}}-1\right)}}.$$
\begin{theorem}
The Nash equilibrium strategies $\beta^*_{j},~j\in\{C,D\}$ related to $\gamma^*=1$ in Table~\ref{table2} and those related to $\gamma^*=0$ in Table~\ref{table3} are acceptable solutions to the APT-I problem with unknown insider threats. 
\end{theorem}

\begin{table*}
	\renewcommand\arraystretch{1.8}
	\caption{\label{table2} Different Nash equilibria in Scenarios C and D}
	\begin{tabular}{|c|c|c|c|c|c|c|}
		\Xhline{1pt} 
		\multicolumn{7}{|c|}{ \textbf{Insider threats are unknown} } \\
		\cline{1-7} 
		\multirow{2}{*}{\textbf{Kinds}} &\multicolumn{3}{c|}{ \textbf{Malicious insider} } &
		\multicolumn{3}{c|}{ \textbf{Inadvertent insider} }   \\
		\cline{2-7} 
		& \multicolumn{2}{c|}{ \textbf{System configurations} } &
		\multirow{1}{*}{\textbf{Nash equilibria}}  & \multicolumn{2}{c|}{ \textbf{System configurations} } &
		\multirow{1}{*}{\textbf{Nash equilibria}} \\
		\cline{1-7} 
		\multirow{2}{*}{\textbf{First}}	& \multirow{2}{*}{$\frac{r_{A}}{p_{A}}=r_{D}\leq1$} & \multirow{2}{*}{$q_{I}\leq (\frac{1}{\frac{r_{A}}{p_{A}}+1})^2-\frac{1}{2}$} &  \multirow{1}{*}{$\left(\frac{r_{D}}{\beta^*_{C}},\beta^*_{C},1\right)$} & \multirow{2}{*}{$r_{A}=r_{D}\leq 1$} & \multirow{2}{*}{$q_{I}\leq (\frac{1}{r_{A}+1})^2-\frac{1}{2}$} & \multirow{1}{*}{$\left(\frac{r_{D}}{\beta^*_{D}},\beta^*_{D},1\right)$}   \\
		&&& with $\beta^*_{C}\in[\varepsilon,1]$ && & with $\beta^*_{D}\in[{\varepsilon},1]$. \\
		\cline{1-7} 
		\textbf{Second}	& $\frac{r_{A}}{p_{A}}\leq1,~\frac{r_{A}}{p_{A}}<r_{D}$ & $q_{I}\leq(\frac{1}{\frac{r_{A}}{p_{A}}+1})^2-\frac{1}{2}$ &  $\left(\frac{r_{A}}{p_{A}},1,1\right)$ & $r_{A}\leq1,~r_{A}<r_{D}$ & $q_{I}\leq (\frac{1}{r_{A}+1})^2-\frac{1}{2}$ & $\left(r_{A},1,1\right)$  \\
		\Xhline{1pt} 
	\end{tabular}
\end{table*}

Table~\ref{table2} shows that $\beta_{C}^*$ are not always identical to $\beta_{D}^*$. It demonstrates that even if the defender does not discover any insiders, the optimal defense strategies for different types of insider threats are different. 

To summarize, the Nash equilibrium strategies $\beta_{j}^*,~j\in\{A,B,C,D\}$ in Tables~\ref{table1} and~\ref{table2} are recommended as acceptable solutions to the APT-I problem with the insider's risky strategy $\gamma^*=1$. Additionally, when the insider avoids risky strategies with $\gamma^*=0$, the strategic interaction between the defender and the attacker is simplified as a two-players APT game. In this case, the Nash equilibria are the same in all Scenarios A-D, and are summarized in Table~\ref{table3}. Thus, the Nash equilibrium strategies $\beta^*$ in Tables~\ref{table3} are recommended as acceptable solutions to the APT-I problem with $\gamma^*=0$, which coincide with the results of the APT problem in~\cite{Hu}.

\begin{table}
	\renewcommand\arraystretch{1.8}
	\caption{\label{table3}The same Nash equilibria in Scenarios A-D}
	\begin{tabular}{|c|c|c|}
		\Xhline{1pt} 
		\multicolumn{2}{|c|}{ \textbf{System configurations} } &
		\multirow{1}{*}{\textbf{Nash equilibria}}    \\
		\cline{1-3}
		\multirow{2}{*}{$r_{A}\!=\!r_{D}\leq1$} & $q_{I}\!>\!(\frac{1}{r_{A}+1})^2\!-\!\frac{1}{2}$ & $\left(\frac{r_{D}}{\beta^*},\beta^*,0\right),~\beta^*\!\in\![r_{D},1]$      \\ 
		\cline{2-3}
		& $q_{I}\!\leq\!(\frac{1}{r_{A}+1})^2\!-\!\frac{1}{2}$ & $\left(\frac{r_{D}}{\beta^*},\beta^*,0\right),~\beta^*\!\in\![r_{D},\varepsilon]$    \\
		\cline{1-3}
		$r_{A}\!\leq\!1,~r_{A}\!<\!r_{D}$ & $q_{I}\!\geq\!(\frac{1}{r_{A}+1})^2\!-\!\frac{1}{2}$  & $(r_{A},1,0)$ \\
		\cline{1-3}
		\multicolumn{2}{|c|}{ $r_{D}\leq1$ and $r_{D}<r_{A}$} & $(1,r_{D},0)$  \\
		\cline{1-3}
		\multicolumn{2}{|c|}{ $r_{A}> 1$ and $r_{D}>1$} & $(1,1,0)$  \\
		\Xhline{1pt}
	\end{tabular}
\end{table}

\begin{remark}
When $\gamma^*=0$, the Nash equilibrium strategies $\alpha^*$ and $\beta^*$ are dependent on $r_{A}$ and $r_{D}$. Specifically, when $r_{A}=r_{D}\leq 1$, $\alpha^*\cdot \beta^* =r_{A}=r_{D}$. When $r_{A}<r_{D}$, the defender focuses on the cost of uncontrolled resources compared with that of the attacker. Hence, the defender reduces the compromised resources by improving his defense strategy. As its counterpart, when $r_{D}<r_{A}$, the attacker infiltrates additional secure resources by improving his attack strategy. When $r_{D}>1$ and $r_{A}>1$, both the defender and the attacker work at full stretch. 
\end{remark}

According to Table~\ref{table3},  $\gamma^*$ is always zero when $r_{D}\leq 1,~r_{D}<r_{A}$ or $r_{A}>1,~r_{D}>1$. In fact, under these conditions, we have $\alpha^*=1$, which lowers the percentage of secure resources $1-x(t)$ on the basis of~\eqref{statesolution}. From~\eqref{costinside}, a small percentage of secure resources causes the profit term $p_{I}(1-x(t))^2+\gamma^2(1-x(t))^2$ to be insufficient to offset the risk term $q_{I}\gamma+\frac{1}{2}\gamma^2$. Thus, under these conditions, the insider is no excuse for action, whose impact can be ignored by the defender. Furthermore, a large risk coefficient $q_{I}>(\frac{1}{r_{A}+1})^2-\frac{1}{2}$ suppresses the risk behaviors of insiders. Based on the above analysis, regardless of the type of insider, the defender can take the following general measures to reduce his costs: i) Reduce the weight of defense action $q_{D}$ so that $r_{D}\geq r_{A}$ or $r_{D}>1$, allowing for large $\beta^*$ and small $\alpha^*$; ii) Use active detection actions to increase $q_{I}$ so that $\gamma^*=0$.

%To sum up, the choices of the insider's optimal strategies are dependent on the percentage of secure resources and the risk coefficient borne by the insider. 

\section{Discussions and Suggestions}\label{Discussion}
\subsection{Discussions}
Regarding two types of insiders in the APT-I game, we compare $\beta^*$ and the corresponding defender's costs. The results show that $\beta^*$ under the malicious threat is not larger than it is under the inadvertent, despite the fact that the malicious insider costs the defender more than the inadvertent.

%Following it, we make some suggestions for the defender when it faces two different types of insider threats. 

\textbf{Optimal defense strategies: }
Table~\ref{table4} lists the corresponding Nash equilibrium strategies of the APT-I game with two types of insiders for the eight system configurations presented in Tables~\ref{table1} and~\ref{table2}.
\begin{table*}
	\centering
	\renewcommand\arraystretch{1.8}
	\caption{\label{table4} The Nash equilibrium strategies under different system configurations and two types of insiders}
	\begin{tabular}{!{\vrule width1pt}c|c|c|c|c!{\vrule width1pt}c|c|c|c|c!{\vrule width1pt}}
		\Xhline{1pt} 
	\textbf{Configurations} & \textbf{Insiders} & $\alpha^*$ & $\beta^*$ & $\gamma^*$ & \textbf{Configurations} & \textbf{Insiders} & $\alpha^*$ & $\beta^*$ & $\gamma^*$\\
		\Xhline{1pt} 
		$\frac{r_{A}}{p_{A}}=\frac{r_{D}}{p_{D}}\leq1~\text{and}$ & Inadvertent & $r_{A}$ & $1$ & $1$ & $\frac{r_{A}}{p_{A}}\leq 1,~\frac{r_{A}}{p_{A}}<\frac{r_{D}}{p_{D}}~\text{and}$ & Inadvertent & $r_{A}$ & $1$ & $1$\\
		\cline{2-5} \cline{7-10} 
		$q_{I}\leq (\frac{1}{\frac{r_{A}}{p_{A}}+1})^2-\frac{1}{2}$	&  Malicious & $[\frac{r_{D}}{p_{D}},\frac{r_{D}}{\varepsilon\sqrt{p_{D}^3}}]~\uparrow$ & $[\varepsilon\sqrt{p_{D}},1]~\downarrow$ & $0$ & $q_{I}\leq (\frac{1}{\frac{r_{A}}{p_{A}}+1})^2-\frac{1}{2}$	&  Malicious & $\frac{r_{A}}{p_{A}}~\uparrow$ & $1$ & $1$ \\
			\cline{1-10} 
		$r_{A}=\frac{r_{D}}{p_{D}}\leq1~\text{and}$ & Inadvertent & $[\frac{r_{D}}{p_{D}},\frac{r_{D}}{\varepsilon\sqrt{p_{D}^3}}]$ & $[\varepsilon\sqrt{p_{D}},1]$ & $1$ & $r_{A}\leq 1,~ r_{A}<\frac{r_{D}}{p_{D}}\leq1~\text{and}$ & Inadvertent & $r_{A}$ & $1$ & $1$\\
		\cline{2-5} \cline{7-10} 
		$q_{I}\leq (\frac{1}{r_{A}+1})^2-\frac{1}{2}$	&  Malicious & $1~\uparrow$ & $r_{D}~\downarrow$ & $0$ & $q_{I}\leq (\frac{1}{r_{A}+1})^2-\frac{1}{2}$	&  Malicious & $\frac{r_{A}}{p_{A}}~\uparrow$ & $1$ & $1$\\
	\cline{1-10} 
	$\frac{r_{A}}{p_{A}}=r_{D}\leq1~\text{and}$ & Inadvertent & $r_{A}$ & $1$ & $1$ & $\frac{r_{A}}{p_{A}}\leq 1,~\frac{r_{A}}{p_{A}}<r_{D}~\text{and}$ & Inadvertent & $r_{A}$ & $1$ & $1$\\
	\cline{2-5} \cline{7-10} 
	$q_{I}\leq (\frac{1}{\frac{r_{A}}{p_{A}}+1})^2-\frac{1}{2}$	&  Malicious & $[r_{D},\frac{r_{D}}{\varepsilon}]~\uparrow$ & $[\varepsilon,1]~\downarrow$ & $1$ & $q_{I}\leq (\frac{1}{\frac{r_{A}}{p_{A}}+1})^2-\frac{1}{2}$	&  Malicious & $\frac{r_{A}}{p_{A}}~\uparrow$ & $1$ & $1$ \\
	\cline{1-10} 
	$r_{A}=r_{D}\leq1~\text{and}$ & Inadvertent & $[r_{D},\frac{r_{D}}{\varepsilon}]$ & $[\varepsilon,1]$ & $1$ & $r_{A}\leq1,~ r_{A}<r_{D}~\text{and}$ & Inadvertent & $r_{A}$ & $1$ & $1$  \\
	\cline{2-5} 	\cline{7-10}
	$q_{I}\leq (\frac{1}{r_{A}+1})^2-\frac{1}{2}$	&  Malicious & $[\frac{r_{D}}{\varepsilon},1]~\uparrow$ & $[r_{D},\varepsilon]~\downarrow$ & $0$ & $q_{I}\leq (\frac{1}{r_{A}+1})^2-\frac{1}{2}$	&  Malicious & $\frac{r_{A}}{p_{A}}~\uparrow$ & $1$ & $1$
	 \\
\Xhline{1pt}
	\end{tabular}
\end{table*}

In comparison, the malicious insider always increases $\alpha^*$ while decreasing $\beta^*$. Hence, $\beta^*$ with an inadvertent insider threat is never less than $\beta^*$ with a malicious insider threat. It suggests that the defender should employ more active defense strategies, such as a higher frequency of user double identifications, network security monitors, and employee training~\cite{applying2,inadvertent3}, to address the APT-I problem with inadvertent insider threats rather than malicious threats.

\textbf{Defender's average costs:}
We compare defender's average costs at Nash equilibria for two types of insider threats.
\subsubsection{Insider threats are known}
In regards to $J_{D}(\alpha^*,\beta^*,\gamma^*)$, the average cost of the defender is calculated by
\begin{flalign}
	J_{D}\left(\alpha^*,\beta^*,1\right)&=(1+q_{D}(\beta^*)^2)(\frac{\alpha^*}{\alpha^*+\beta^*})^{2},\nonumber\\
	&=\frac{1}{1+q_{D}(\beta^*)^2},~\label{j1}
\end{flalign}
since $\alpha^*\cdot \beta^*=\frac{r_{D}}{p_{D}}$ and $\gamma^*=1$ according to Table~\ref{table1}.

\subsubsection{Insider threats are unknown}
The average cost of the defender is calculated by
\begin{flalign}
	J_{D}\left(\alpha^*,\beta^*,0\right)&=(p_{D}+q_{D}(\beta^*)^2)(\frac{\alpha^*}{\alpha^*+\beta^*})^{2},\nonumber\\
	&=\frac{p_{D}^2}{p_{D}+q_{D}(\beta^*)^2},~\label{j2}
\end{flalign}
since $\alpha^*\cdot \beta^*=r_{D}$ according to Table~\ref{table2}.

Based on~\eqref{j1} and~\eqref{j2}, the malicious insider threat causes more harm to the defender than the inadvertent insider threat, because the malicious insider always reduces $\beta^*$, as shown in Table~\ref{table4}.

\subsection{Suggestions}~\label{ObservationD}
The Nash equilibrium strategies $\beta^*$ in Tables~\ref{table1}-\ref{table3} are recommended as acceptable solutions to the APT-I problem. Furthermore, the following suggestions can be used by the defender to achieve a much more feasible defense against the APT with two types of insider threats.
\begin{itemize}
\item In the APT-I problem, $\beta^*$ with an inadvertent insider threat is never less than $\beta^*$ with a malicious insider threat, which suggests that the defender should employ more active defense strategies for inadvertent insider threats than for malicious threats. Furthermore, if the defender knows noting about insider threats, he can employ a relatively high level of monitoring and defense.

\item A sufficiently large $q_{I}>(\frac{1}{\frac{r_{A}}{p_{A}}+1})^2-\frac{1}{2}$ can always ensure $\gamma^*=0$. As a result, the defender can increase the intensity of regulation and the frequency of monitoring to raise $q_{I}$ such that $\gamma^*=0$. For the inadvertent insider, the defender can increase employee training and risk communication, engender process discipline, and heighten punitive measures to stop the inadvertent insider's risky behaviors~\cite{inadvertent3}. For the malicious insider, the defender can use deception techniques, such as Darknets, Honey nets, to increase the risk of being detected when the malicious insider copies and sells sensitive information~\cite{detection}.

\item When $r_{D}\leq 1,~r_{D}<r_{A}$ or $r_{D}>1,~r_{A}>1$, $\gamma^*$ is always zero. It implies that under those conditions, the defender can ignore the effect of insider threats and maintain the same defense strategy and monitoring intensity as in the APT game.
\end{itemize}

\section{Numerical Experiments}\label{Example}
In this section, we first check the Nash equilibria in Table~\ref{table1}-\ref{table3}. Then we show the discussion presented in Section~\ref{Discussion}. Finally, by examining the conditions of risky strategies adopted by different insiders, we demonstrate that the inadvertent insider is more likely to lead to the risky strategy.

\subsection{The Nash equilibrium verification}
To validate the Nash equilibria, we investigate the minimums of the average cost functions of the defender and the attacker, as well as the maximums of the average profit functions of insiders. 

\subsubsection{Verification of Tables~\ref{table1} and~\ref{table2}}
Let $p_{A}=0.95$, $p_{D}=0.9$, $q_{D}=4.5$, $p_{I}=0.1$ and $q_{I}=0.01$. The setting of $q_{A}$ and Nash equilibria are listed in Table~\ref{table5}. The results are shown in Fig. 2 and Fig. 3.
	\begin{table}\tiny
	\centering
	\renewcommand\arraystretch{1.8}
	\caption{\label{table5}The system configurations of Tables~\ref{table1} and~\ref{table2}}
	\begin{tabular}{c|c|c|c|c}
		\Xhline{1pt} 
		\textbf{Scenarios} & \textbf{Kinds} & $q_{A}$ & \textbf{System configurations} & \textbf{Nash Equilibria}  \\
		\Xhline{1pt}
		\multirow{2}{*}{A} & First & $4.55$ & $\frac{r_{A}}{p_{A}}\!=\!\frac{r_{D}}{p_{D}}\!=\!0.22 \!\leq\! 1$ & $(0.26,0.84,1)$      \\ 
    	\cline{2-5}
	   & Second & $7.14$ & $\frac{r_{A}}{p_{A}}\!=\!0.14\!<\!\frac{r_{D}}{p_{D}}\!=\!0.22 $ & $(0.14,1,1)$\\
	   	\cline{1-5}
	   \multirow{2}{*}{B} & First & $4.32$ & $r_{A}\!=\!\frac{r_{D}}{p_{D}}\!=\!0.22 \!\leq\! 1$ & $(0.25,0.89,1)$      \\ 
	   \cline{2-5}
	   & Second & $7.14$ & $r_{A}\!=\!0.13\!<\!\frac{r_{D}}{p_{D}}\!=\!0.22 $ & $(0.13,1,1)$ \\
		\cline{1-5}
			\multirow{2}{*}{C} & First & $5$ & $\frac{r_{A}}{p_{A}}\!=\!r_{D}\!=\!0.2\! \leq\! 1$ & $(0.26,0.76,1)$      \\ 
		\cline{2-5}
		& Second & $6.67$ & $\frac{r_{A}}{p_{A}}\!=\!0.15\!<\!r_{D}\!=\!0.2 $ & $(0.15,1,1)$\\
	    \cline{1-5}
	    	\multirow{2}{*}{D} & First & $4.75$ & $r_{A}\!=\!r_{D}\!=\!0.2 \leq 1$ & $(0.25,0.8,1)$      \\ 
	    \cline{2-5}
	    & Second & $6.67$ & $r_{A}\!=\!0.14\!<\!r_{D}\!=\!0.2 $ & $(0.14,1,1)$\\
	
		\Xhline{1pt}
	\end{tabular}
\end{table}

\begin{figure*}
	\centering
	\subfigure[The APT-I game with a known malicious insider (Scenario A)]{\includegraphics[height=4cm,width=15cm]{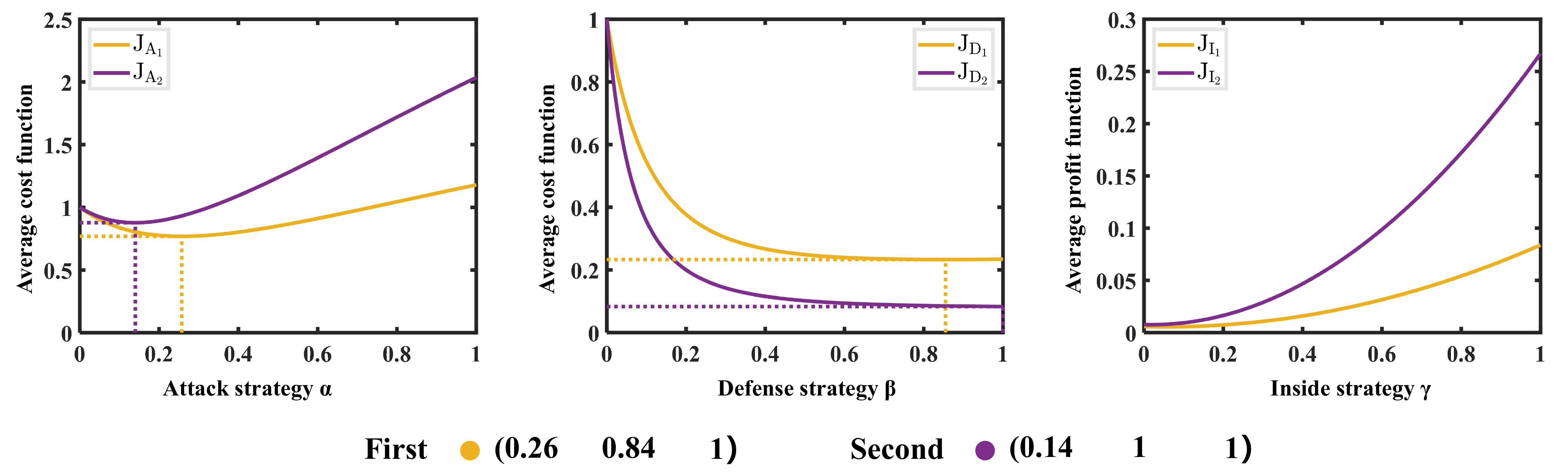}}\\
	\subfigure[The APT-I game with a known inadvertent insider (Scenario B)]{\includegraphics[height=4cm,width=15cm]{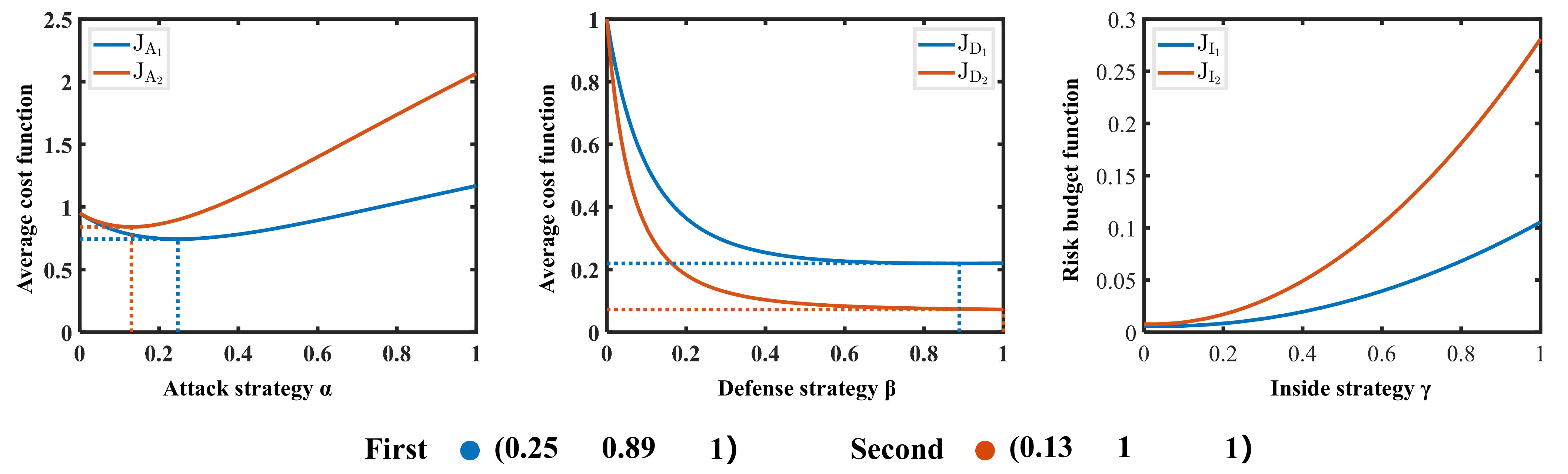}}
	\caption{Verification of Table~\ref{table1}. In Fig. 2(a), the yellow lines show that the minimums of $J_{D}$ and $J_{A}$, respectively, and the maximum of $J_{I}$ fall into the Nash equilibrium $(0.26,0.84,1)$, and the purple lines show that those fall into $(0.14,1,1)$. In Fig. 2(b), the blue lines show that the minimums of $J_{D}$ and $J_{A}$, respectively, and the maximum of $J_{I}$ fall into the Nash equilibrium $(0.25,0.89,1)$, and the red lines show that those fall into $(0.13,1,1)$.}
\end{figure*}

\begin{figure*}
	\centering
	\subfigure[The APT-I game with an unknown malicious insider (Scenario C)]{\includegraphics[height=5cm,width=8.4cm]{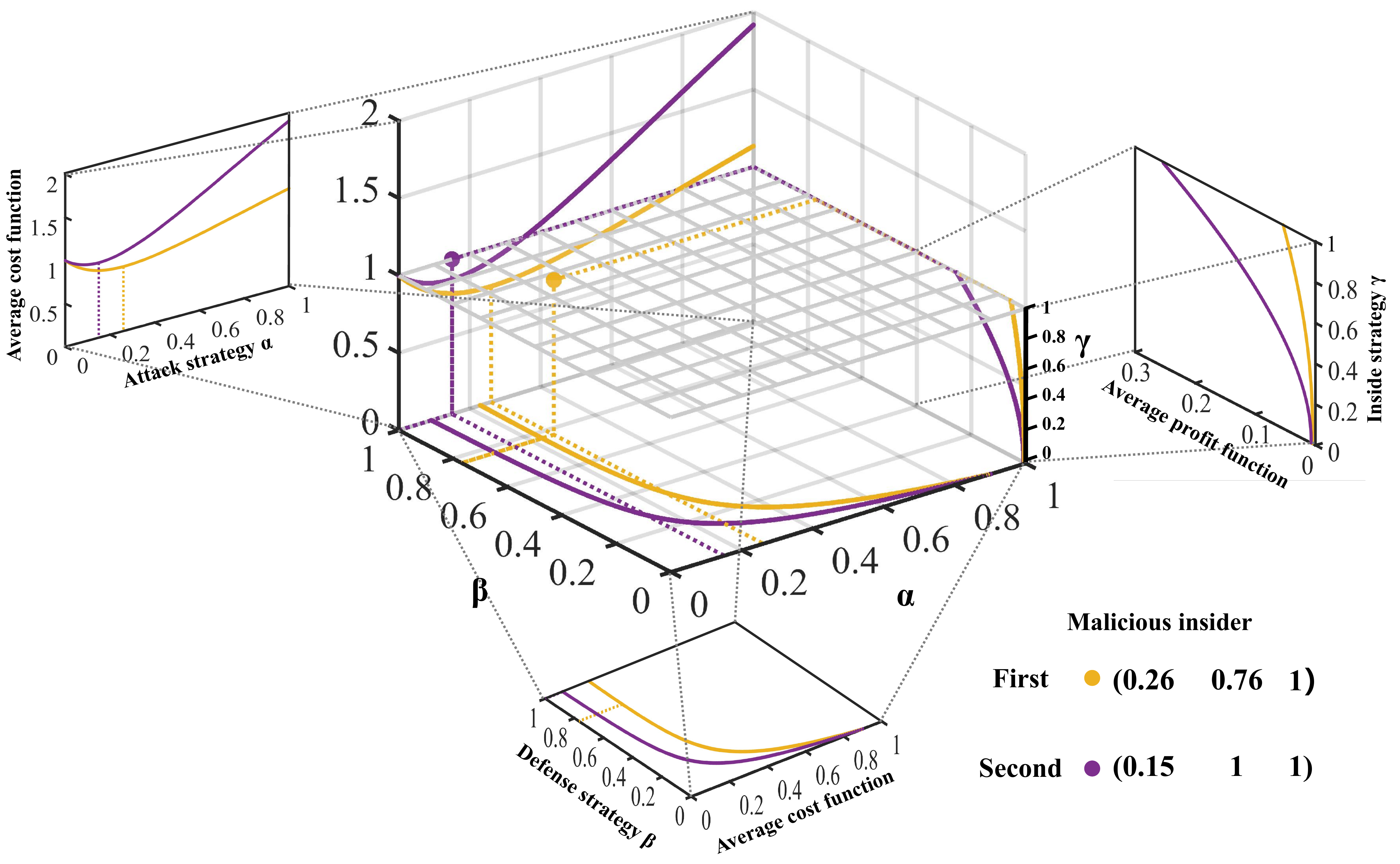}}~~~~~~~~~~~~~~
	\subfigure[The APT-I game with an uknown malicious insider (Scenario D)]{\includegraphics[height=5cm,width=8.4cm]{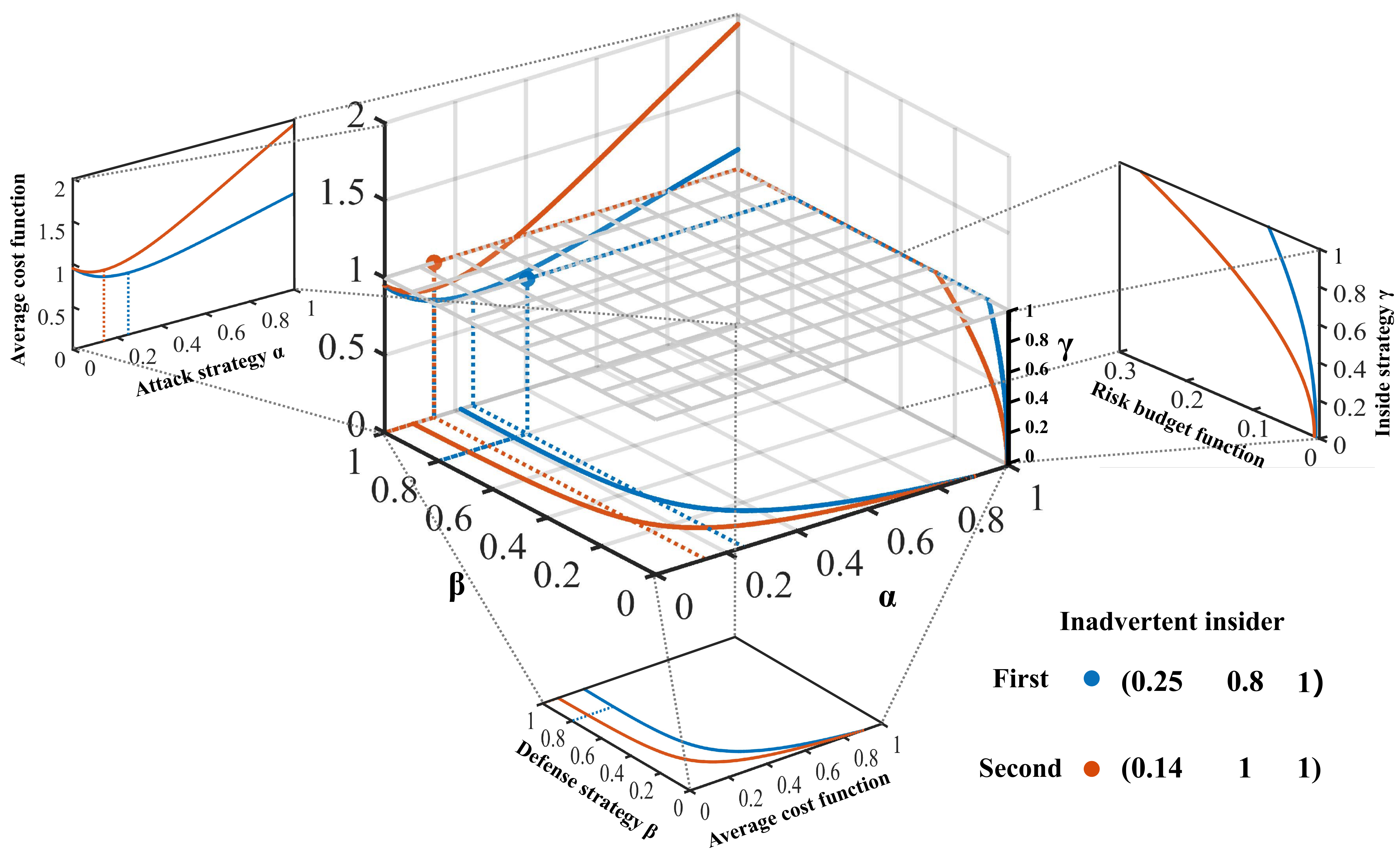}}\\
	\caption{Verification of Table~\ref{table2}. In Fig. 3(a), the yellow lines show that the minimums of $J_{D}$ and $J_{A}$, respectively, and the maximum of $J_{I}$ fall into the Nash equilibrium $(0.26,0.76,1)$ and the purple lines show that those fall into $(0.15,1,1)$. In Fig. 3(b), the blue lines show that the minimums of $J_{D}$ and $J_{A}$, respectively, and the maximum of $J_{I}$ fall into the Nash equilibrium $(0.25,0.8,1)$ and the red lines show that those fall into $(0.14,1,1)$.  }
\end{figure*}

Fig. 2 depicts the average objective function of each player under the system configurations in Table~\ref{table1}. In Fig. 2(a), the Nash equilibrium $(0.26,0.84,1)$ meets $(\frac{r_{D}}{p_{D}\beta^*},\beta^*,1)$ with $\beta^*\in[\varepsilon\sqrt{p_{D}},1]$ and $\varepsilon\sqrt{p_{D}}=0.67$, and the Nash equilibrium $(0.14,1,1)$ meets $(\frac{r_{A}}{p_{A}},1,1)$.  Similarly, in Fig. 2(b), the Nash equilibrium $(0.25,0.89,1)$ satisfies $(\frac{r_{D}}{p_{D}\beta^*},\beta^*,1)$ with $\beta^*\in[\varepsilon\sqrt{p_{D}},1]$ and $\varepsilon\sqrt{p_{D}}=0.67$, and the Nash equilibrium $(0.13,1,1)$ satisfies $(r_{A},1,1)$. 

Fig. 3 draws the average objective function of each player under the system configurations in Table~\ref{table2}. In Fig. 3(a), the Nash equilibrium $(0.26,0.76,1)$ meets $(\frac{r_{D}}{\beta^*},\beta^*,1)$ with $\beta^*\in[\varepsilon,1]$ and $\varepsilon=0.71$, and the Nash equilibrium $(0.15,1,1)$ meets $(\frac{r_{A}}{p_{A}},1,1)$. Similarly, in Fig. 3(b), the Nash equilibrium $(0.25,0.8,1)$ satisfies $(\frac{r_{D}}{\beta^*},\beta^*,1)$ with $\beta^*\in[\varepsilon,1]$ and $\varepsilon=0.71$, and the Nash equilibrium $(0.14,1,1)$ satisfies $(r_{A},1,1)$.

Therefore, the minimum of $J_{D}$, the minimum of $J_{A}$, and the maximum of $J_{I}$ in Fig. 2 and Fig. 3, correspond to each Nash equilibrium in Table~\ref{table5}. 

\subsubsection{Verification of Table~\ref{table3}}\label{Same}
Let $p_{A}=0.95$, $p_{D}=0.9$ and $p_{I}=0.1$. The rest of system configurations and Nash equilibria are set in Table~\ref{table6}. The verification of Nash equilibria is shown in Fig. 4, where the insider always benefits the most when $\gamma^*=0$.

\begin{table}\tiny
	\centering
	\renewcommand\arraystretch{1.8}
	\caption{\label{table6}The system configurations of Tables~\ref{table3}}
	\begin{tabular}{c|c|c|c|c}
		\Xhline{1pt} 
	      $q_{A}$ & $q_{D}$ &	\multicolumn{2}{|c|}{\textbf{System configurations}} & \textbf{Nash Equilibria}   \\
		\Xhline{1pt}
		\multirow{2}{*}{$4.75$} & \multirow{2}{*}{$4.5$} & \multirow{2}{*}{$r_{A}=r_{D}=0.2\leq 1$} & $q_{I}=0.3$ &  $(0.9,0.33,0)$   \\ 
		\cline{4-5}
		 &  & & $q_{I}=0.01$ &  $(0.6,0.5,0)$   \\ 
    	\cline{1-5}
		$6.67$ & $4.5$ & $r_{A}=0.14<r_{D}=0.2$ & $q_{I}=0.3$ &  $(0.225,1,0)$   \\ 
		\cline{1-5}
		$1.9$ & $4.5$ & $r_{D}=0.2<r_{A}=0.5$ & $q_{I}=0.01$ &  $(1,0.3,0)$   \\ 
	  \cline{1-5}
			$0.5$ & $0.3$ & $r_{A}=1.9,~r_{D}=3$ & $q_{I}=0.01$ &  $(1,1,0)$   \\ 
		\Xhline{1pt}
	\end{tabular}
\end{table}

\begin{figure}
	\centering
	\includegraphics[height=7cm,width=8.7cm]{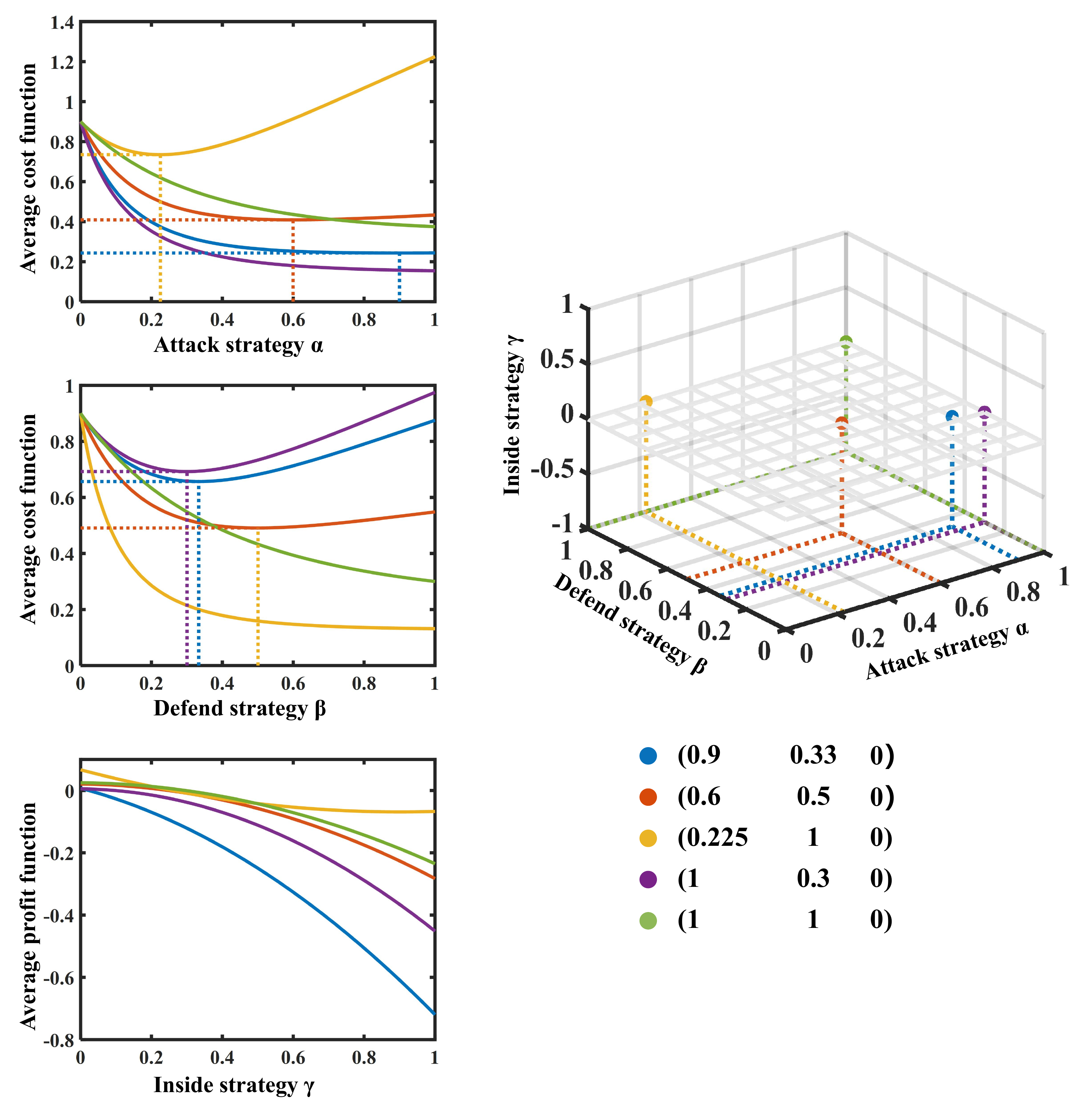}
	\caption{Verification of Table~\ref{table3} with $\gamma^*=0$.  }
\end{figure}

\subsection{Discussion examinations}
\subsubsection{Nash equilibrium strategies}
The Nash equilibrium strategies of the APT-I game with two types of insiders are shown in Table~\ref{table7} under the system configurations in Table~\ref{table5}. It demonstrates that $\beta^*$ related to the malicious threat is never larger than $\beta^*$ related to the inadvertent threat.
\begin{table*}\tiny
	\centering
	\renewcommand\arraystretch{1.8}
	\caption{\label{table7} The variations of Nash equilibrium strategies}
	\begin{tabular}{c|c|c|c|c|c|c|c|c|c}
		\Xhline{1pt} 
		\textbf{System configurations} & \textbf{Insider Types} & $\alpha^*$ & $\beta^*$ & $\gamma^*$ & \textbf{System configurations} & \textbf{Insider Types} & $\alpha^*$ & $\beta^*$ & $\gamma^*$\\
		\Xhline{1pt} 
		$\frac{r_{A}}{p_{A}}=\frac{r_{D}}{p_{D}}=0.22\leq1~\text{and}$ & Inadvertent & $0.21$ & $1$ & $1$ & $\frac{r_{A}}{p_{A}}=0.14<\frac{r_{D}}{p_{D}}=0.22~\text{and}$ & Inadvertent & $0.13$ & $1$ & $1$\\
		\cline{2-5}\cline{7-10} 
		$q_{I}=0.01< 0.17$	&  Malicious & $[0.22,0.33]~\uparrow$ & $[0.67,1]~\downarrow$ & $0$ & $q_{I}=0.01< 0.27$	&  Malicious & $0.14~\uparrow$ & $1$ & $1$ \\	
			\cline{1-10} 
		$r_{A}=\frac{r_{D}}{p_{D}}=0.22\leq1~\text{and}$ & Inadvertent & $[0.22,0.33]$ & $[0.67,1]$ & $1$ & $r_{A}=0.13<\frac{r_{D}}{p_{D}}=0.22~\text{and}$ & Inadvertent & $0.13$ & $1$ & $1$\\
		\cline{2-5}\cline{7-10} 
		$q_{I}=0.01<0.17$	&  Malicious & $1~\uparrow$ & $0.2~\downarrow$ & $0$ & $q_{I}=0.01<0.28$	&  Malicious & $0.14~\uparrow$ & $1$ & $1$\\
    	\cline{1-10} 
		$\frac{r_{A}}{p_{A}}=r_{D}=0.2\leq1~\text{and}$ & Inadvertent & $0.19$ & $1$ & $1$ & $\frac{r_{A}}{p_{A}}=0.15<r_{D}=0.2~\text{and}$ & Inadvertent & $0.14$ & $1$ & $1$\\
		\cline{2-5}\cline{7-10} 
		$q_{I}=0.01<0.19$	&  Malicious & $[0.2,0.28]~\uparrow$ & $[0.71,1]~\downarrow$ & $1$ & $q_{I}=0.01<0.26$	&  Malicious & $0.15~\uparrow$ & $1$ & $1$ \\
			\cline{1-10} 
	$r_{A}=r_{D}=0.2\leq1~\text{and}$ & Inadvertent & $[0.2,0.28]$ & $[0.71,1]$ & $1$ & $ r_{A}=0.14<r_{D}=0.2~\text{and}$ & Inadvertent & $0.14$ & $1$ & $1$  \\
	\cline{2-5}\cline{7-10} 
	$q_{I}=0.01< 0.19$	&  Malicious & $[0.28,1]~\uparrow$ & $[0.2,0.71]~\downarrow$ & $0$ & $q_{I}=0.01<0.27$	&  Malicious & $0.15~\uparrow$ & $1$ & $1$ \\
		\Xhline{1pt} 
	\end{tabular}
\end{table*}

\subsubsection{Average costs of all players}
In Fig. 5, we compare the defender's average costs with malicious insider threats to those with inadvertent insider threats. The results show that the malicious insider costs the defender more than the inadvertent insider. Similar comparisons for the attacker's average cost and the insider's average profit are shown in Fig. 6. According to the findings, the malicious insider always raises the attacker's cost, whose profit is lower than the inadvertent insider's profit.
\begin{figure}
	\centering
	\subfigure[Comparisons of $J_{D}$ between two types of known insider threats ($J_{D}$ in Scenarios A and B)]{\includegraphics[height=3cm,width=4.3cm]{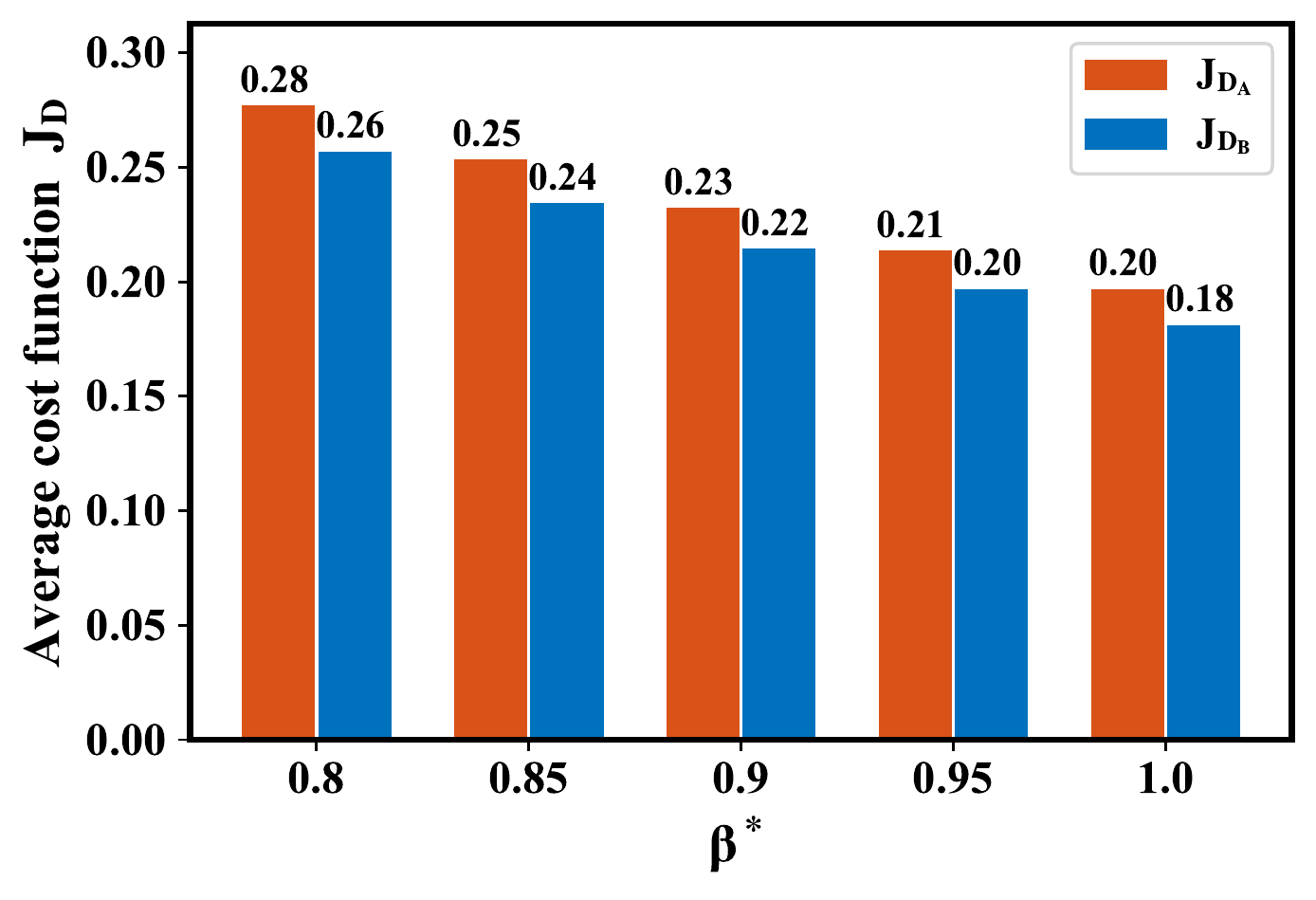}}
	\subfigure[Comparisons of $J_{D}$ between two types of unknown insider threats ($J_{D}$ in Scenarios C and D)]{\includegraphics[height=3cm,width=4.3cm]{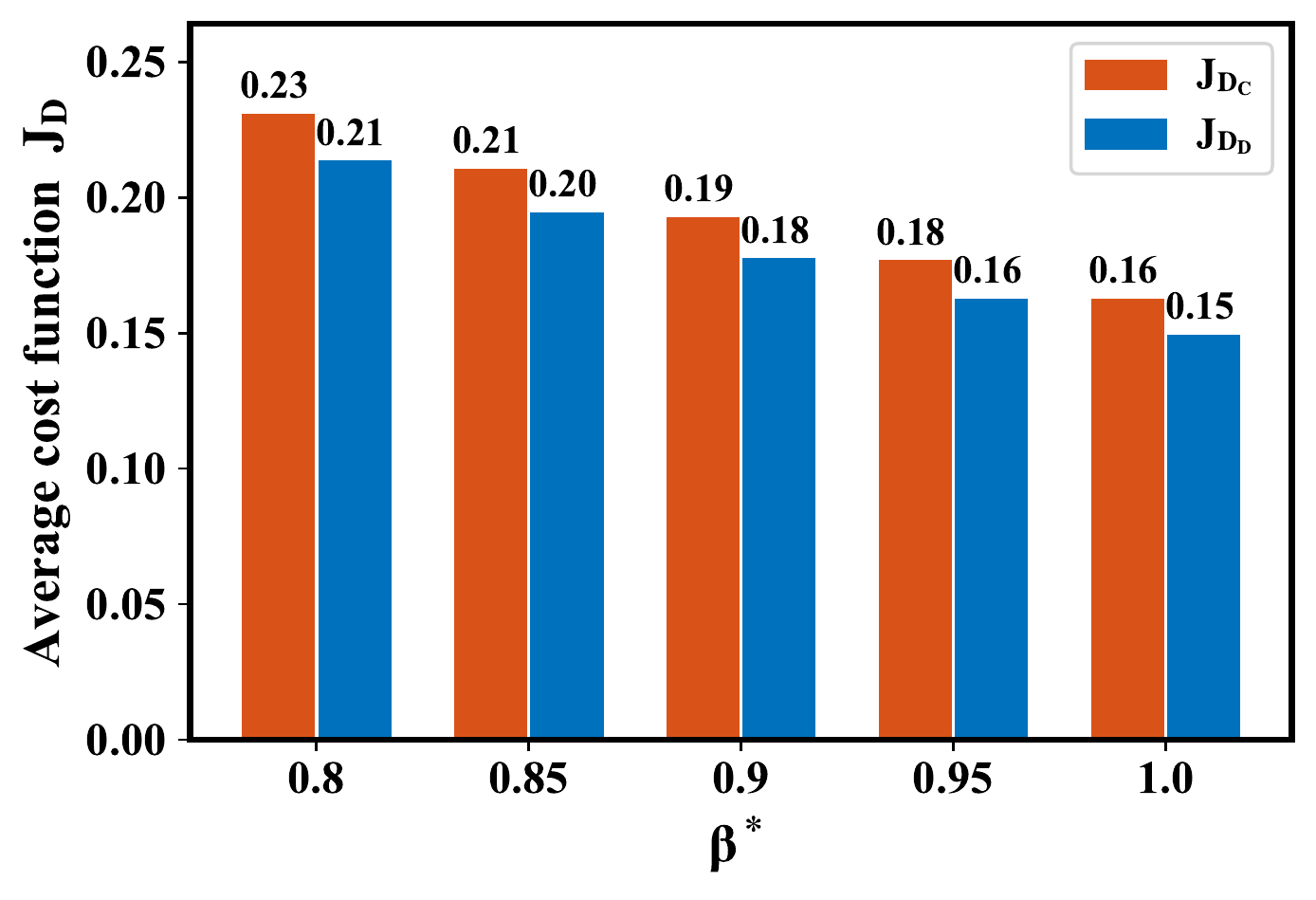}}
	\caption{Comparisons of the defender's average costs under two types of insider threats. The red bars represent the values of $J_{D}$ under various $\beta^*$ associated with a malicious insider, and the blue bars represent those associated with an inadvertent insider. }
\end{figure}

\begin{figure}
	\centering
	\subfigure[Comparisons of $J_{A}$ and $J_{I}$ between two types of known insider threats, respectively ($J_{A}$ and $J_{I}$ in Scenarios A and B)]{\includegraphics[height=3cm,width=4.3cm]{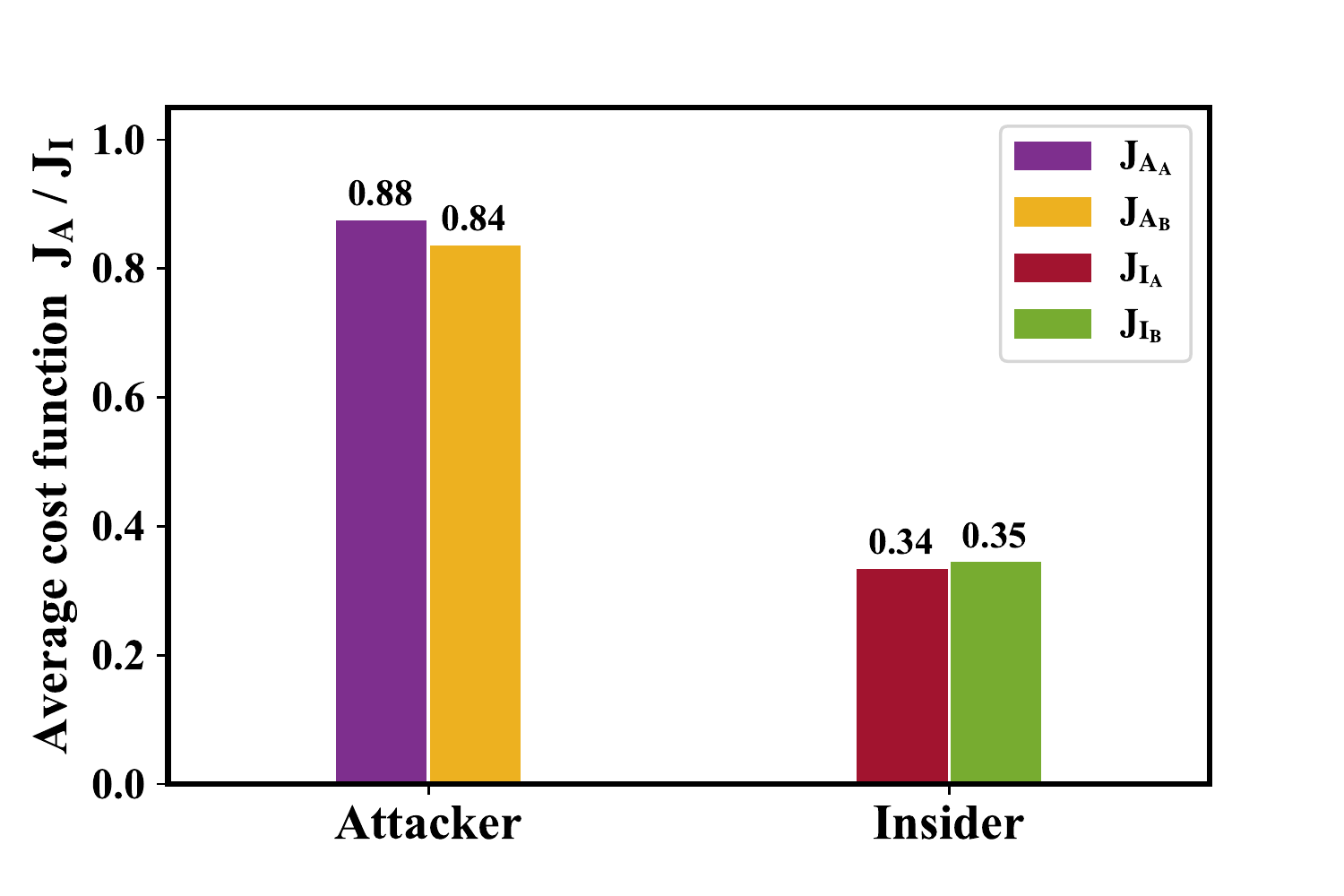}}\label{Igamma1I.1}
	\subfigure[Comparisons of $J_{A}$ and $J_{I}$ between two types of unknown insider threats, respectively ($J_{A}$ and $J_{I}$ in Scenarios C and D)]{\includegraphics[height=3cm,width=4.3cm]{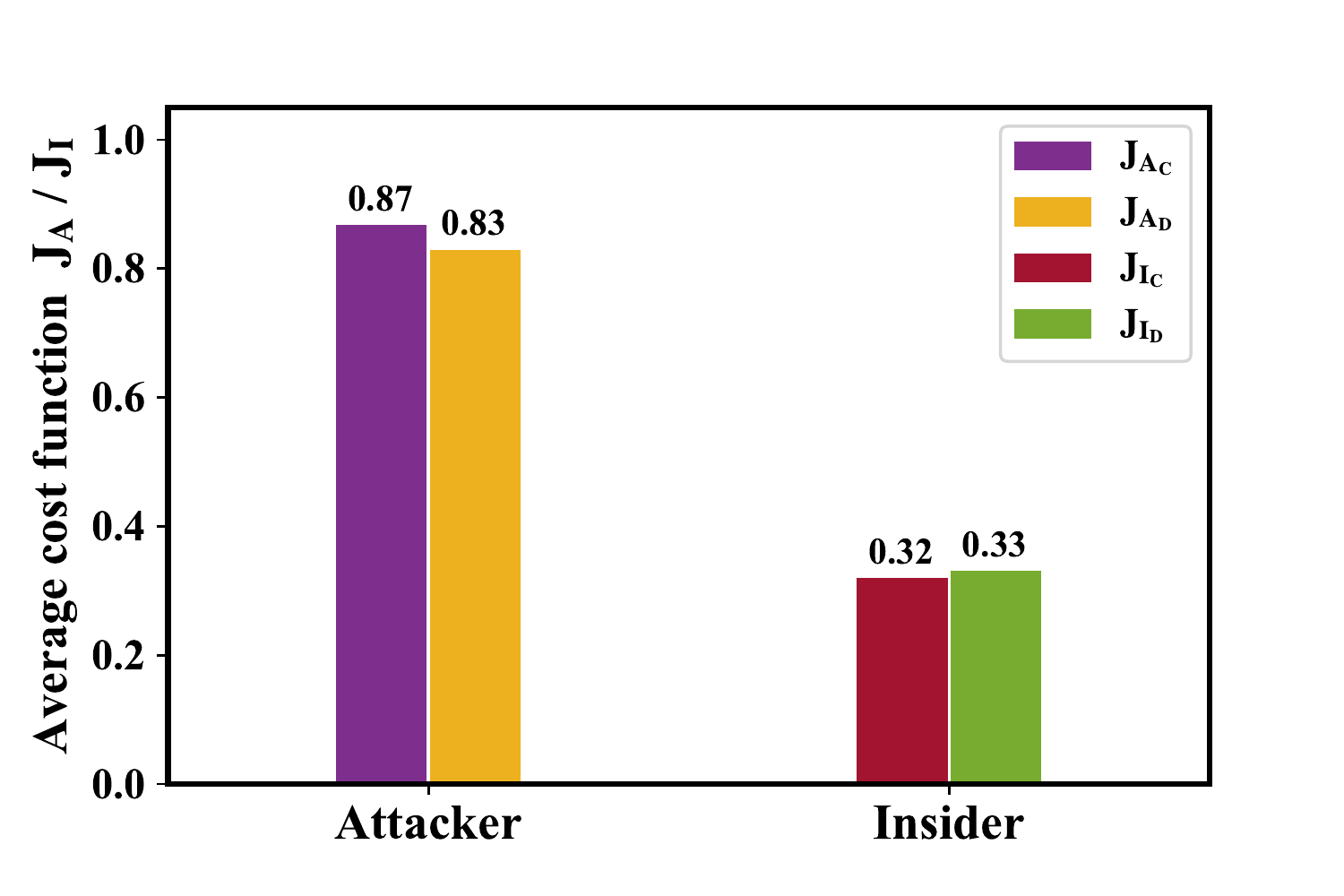}}\label{Igamma1I.2}\\
	\caption{Comparisons of the attacker's average costs and the insider's average profits under two types of insider threats. The purple bars represent the values of $J_{A}$ associated with a malicious insider, and the yellow bars represent those associated with an inadvertent insider. The red bars represent the values of $J_{I}$ associated with a malicious insider, and the green bars represent those associated with an inadvertent insider.}
\end{figure}

\subsection{The examination of  risk coefficients}\label{B33}
Different risk coefficients $q_{I}$ result in different Nash equilibrium strategies $\gamma^*$ for the insider. Hence, we set a series of $q_{I}^{k}=0.1+0.01*k,~k=0,\cdots,14$ to examine the conditions of the risky strategy pursued by different insiders. Let $p_{A}=p_{D}=0.8$, $q_{A}=4$ and $p_{I}=0.1$. The rest of the parameters and available Nash equilibria are listed in Table~\ref{table8}. According to Fig. 7, the inadvertent insider is more likely than the malicious insider to maximize profit at $\gamma^*=1$, implying that the inadvertent insider is more likely to lead to a risky strategy.
As a result, the defender must employ more monitoring efforts (a larger $q_{I}$) to avoid the risky behavior of the inadvertent insider compared with the malicious insider.
	\begin{table}\tiny
	\centering
	\renewcommand\arraystretch{1.8}
	\caption{\label{table8}The system configurations with risk coefficients}
	\begin{tabular}{c|c|c|c|c}
		\Xhline{1pt} 
		\textbf{Scenarios} & $q_{D}$ & $r_{D}$ & $q_{I}$ & \textbf{Nash Equilibria}  \\
		\Xhline{1pt}
			\multirow{3}{*}{A} & \multirow{3}{*}{$4$} & \multirow{3}{*}{$0.2$} & $q_{I}\leq 0.14$ & $(0.25,1,1)$      \\ 
		\cline{4-5}
		& & & $0.14<q_{I}<0.19$   &$(1,0.2,0)$\\
		\cline{4-5}
		& & & $q_{I}\geq0.19$ & $(0.2,1,0)$\\
		\cline{1-5}
	\multirow{2}{*}{B} & \multirow{2}{*}{$5$} & \multirow{2}{*}{$0.16$}& $q_{I}\leq 0.19$ & $(0.2,1,1)$      \\ 
	\cline{4-5}
	& & & $q_{I}>0.19$ & $(1,0.16,0)$\\
		\cline{1-5}
		\multirow{3}{*}{C} & \multirow{3}{*}{$3.2$} & \multirow{3}{*}{$0.25$}& $q_{I}\leq 0.14$ & $(0.25,1,1)$      \\ 
		\cline{4-5}
		& & & $0.14<q_{I}<0.19$ & None\\
		\cline{4-5}
		& & & $q_{I}\geq 0.19$ & $(0.2,1,0)$\\
			\cline{1-5}
		\multirow{2}{*}{D} & \multirow{2}{*}{$4$} &  \multirow{2}{*}{$0.2$}& $q_{I}\leq 0.19$ & $(0.2,1,1)$      \\ 
	\cline{4-5}
	& & &$q_{I}>0.19$ & $(0.2,1,0)$\\
		\Xhline{1pt}
	\end{tabular}
\end{table}

\begin{figure}
	\centering
	\subfigure[$J_{I}(1)-J_{I}(0)$ in Scenario A and Scenario B]{\includegraphics[height=3cm,width=4.3cm]{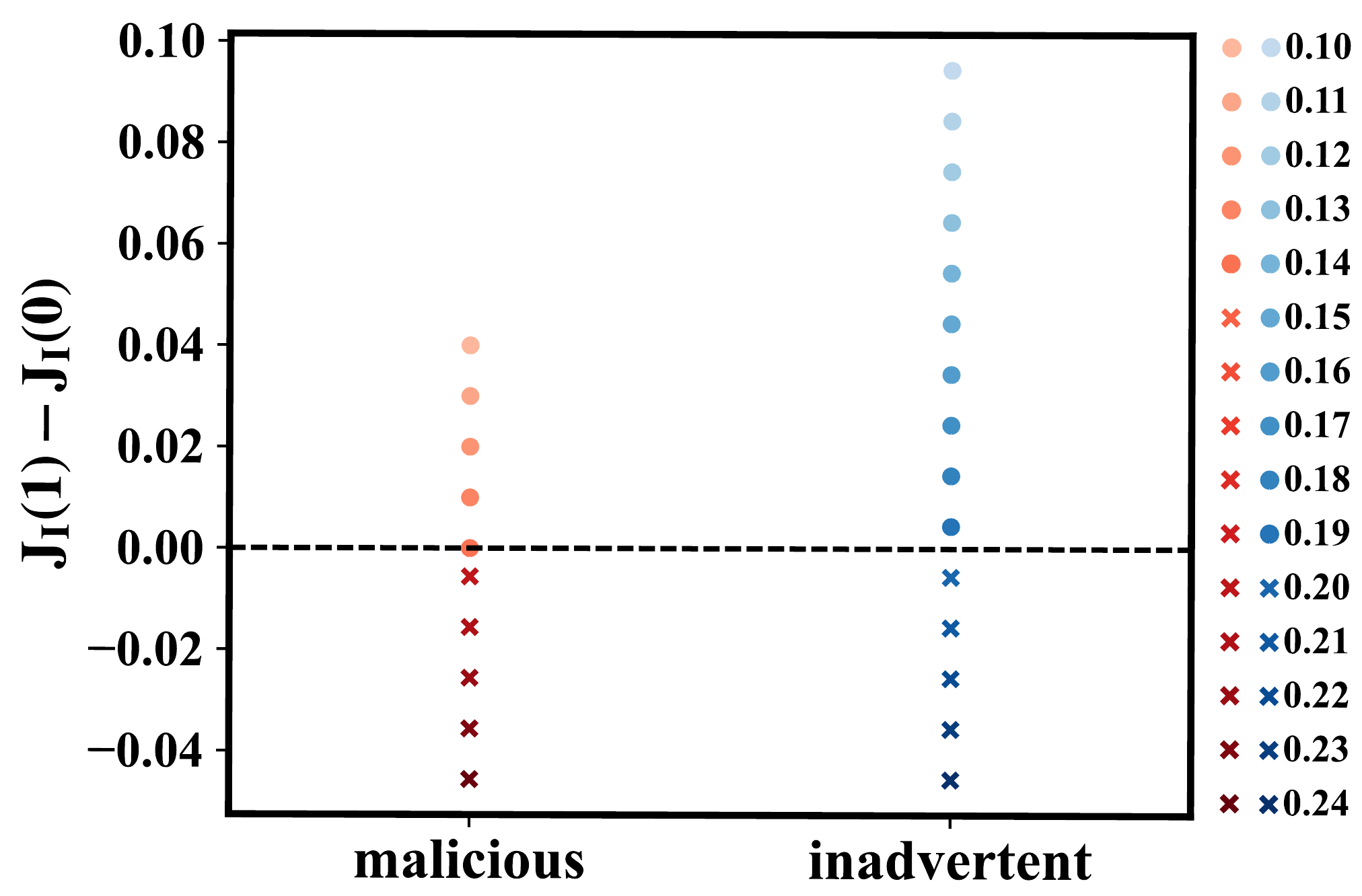}}
	\subfigure[$J_{I}(1)-J_{I}(0)$ in Scenario C and Scenario D]{\includegraphics[height=3cm,width=4.3cm]{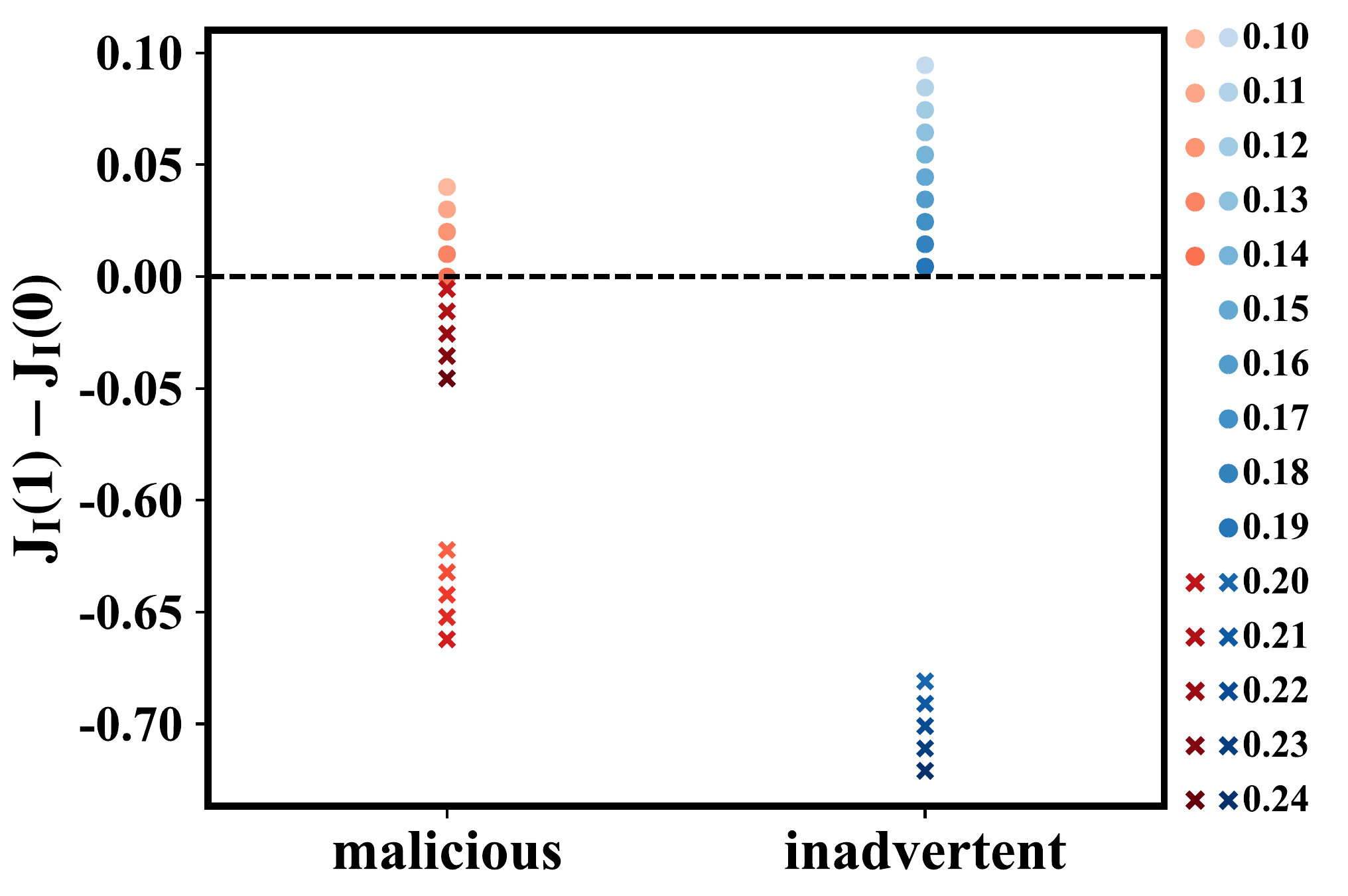}}\\
	\caption{The influence of different $q_{I}$ for Nash equilibrium strategies $\gamma^*$. If $J_{I}(1)-J_{I}(0)\geq 0$, then $\gamma^*=1$, otherwise $\gamma^*=0$.}
\end{figure}

\section{Conclusion}\label{Conclusion}
This paper investigated the APT-I problem in the context of malicious insider threats or inadvertent insider threats. To address this problem, we developed a unified three-player APT-I game framework. To deal with the two types of insider threats, we obtained feasible Nash equilibria, in which the optimal defense strategies were recommended as acceptable solutions to the APT-I problem. Furthermore,  the optimal defense strategies and the average cost of the defender have been compared between two types of insider threats. The results demonstrated the differences between two types of insider threats in the APT-I problem, and we further provided a list of effective suggestions for the defender.
\appendices
\section{Proof of Lemma 1}\label{proofbest}
By plugging \eqref{statesolution} into \eqref{costattack}-\eqref{costinside}, we obtain
\begin{flalign}
&J_{A}\left(\alpha,\beta,\gamma\right)= (p_{A}(1-\gamma)^2+q_{A}\alpha^2+\gamma^2)(\frac{\beta}{\alpha+\beta})^2,\nonumber\\
&J_{D}\left(\alpha,\beta,\gamma\right)=(p_{D}(1-\gamma)^2+q_{D}\beta^2+\gamma^2)(\frac{\alpha}{\alpha+\beta})^{2},\nonumber\\
&J_{I}\left(\alpha,\beta,\gamma\right)= (p_{I}+\gamma^2)(\frac{\beta}{\alpha+\beta})^{2}-(q_{I}\gamma+\frac{1}{2}\gamma^2).\nonumber
\end{flalign}
Taking the derivative of $J_{A}(\alpha,\beta^*,\gamma^*)$, we obtain
\begin{eqnarray}
	\frac{d{J}_{A}(\alpha,\beta^*,\gamma^*)}{d\alpha}=\frac{2(\beta^*)^2}{(\alpha+\beta^*)^3}[q_{A}\alpha\beta^*-p_{A}(1-\gamma^*)^2-(\gamma^*)^2].\nonumber
\end{eqnarray}
When $\alpha=0$, it follows that
\begin{eqnarray}
\frac{d{J}_{A}(\alpha,\beta^*,\gamma^*)}{d\alpha}\Big{|}_{\alpha=0}=-\frac{2}{\beta^*}[p_{A}(1-\gamma^*)^2+(\gamma^*)^2]<0,\nonumber
\end{eqnarray}
and therefore, $J_{A}(\alpha,\beta^*,\gamma^*)$ is decreasing at the initial attack strategy. With $\alpha\in (0,1]$, if $\frac{d{J}_{A}(\alpha,\beta^*,\gamma^*)}{d\alpha}$ stays negative, ${J}_{A}(\alpha,\beta^*,\gamma^*)$ reaches the minimum value at $\alpha^*=1$. If $\frac{d{J}_{A}}{d\alpha}$ strikes the horizontal axis once and only once at the point $\bar{\alpha}=\frac{p_{A}(1-\gamma^*)^2+(\gamma^*)^2}{q_{A}\beta^*}$, then the minimum of ${J}_{A}(\alpha,\beta^*,\gamma^*)$ is $\alpha^*=\bar{\alpha}$ with $\bar{\alpha}<1$ or $\alpha^*=1$ with $\bar{\alpha}\geq 1$.

Similarly, taking the derivative of $J_{D}(\alpha^*,\beta,\gamma^*)$, we obtain
\begin{eqnarray}
	\frac{d{J}_{D}(\alpha^*,\beta,\gamma^*)}{d\beta}=\frac{2(\alpha^*)^2}{(\alpha^*+\beta)^3}[q_{D}\alpha^*\beta-p_{D}(1-\gamma^*)^2-(\gamma^*)^2].\nonumber
\end{eqnarray}
When $\beta=0$, it follows that
\begin{eqnarray}
	\frac{d{J}_{D}(\alpha^*,\beta,\gamma^*)}{d\beta}\Big{|}_{\beta=0}=-\frac{2}{\alpha^*}[p_{D}(1-\gamma^*)^2+(\gamma^*)^2]<0.\nonumber
\end{eqnarray}
Thus, if $\bar{\beta}=\frac{p_{D}(1-\gamma^*)^2+(\gamma^*)^2}{q_{D}\alpha^*}<1$, then the minimum of $J_{D}(\alpha^*,\beta,\gamma^*)$ is $\beta^*=\bar{\beta}$, otherwise $\beta^*=1$.

Finally, we derive the best response strategy of the insider. Taking the derivative of ${J}_{I}(\alpha^*,\beta^*,\gamma)$, we obtain
\begin{eqnarray}
	\frac{d{J}_{I}(\alpha^*,\beta^*,\gamma)}{d\gamma}=2\gamma(\frac{\beta^*}{\alpha^*+\beta^*})^2-q_{I}-\gamma.\nonumber
\end{eqnarray}
Furthermore, the second derivative of ${J}_{I}(\alpha^*,\beta^*,\gamma)$ satisfies
\begin{eqnarray} \frac{d^2{J}_{I}(\alpha^*,\beta^*,\gamma)}{d^2\gamma}=2(\frac{\beta^*}{\alpha^*+\beta^*})^2-1,\nonumber
\end{eqnarray}	
which indicates that if $2(\frac{\beta^*}{\alpha^*+\beta^*})^2-1\geq0$, ${J}_{I}(\alpha^*,\beta^*,\gamma)$ is convex, otherwise it would be concave. Meanwhile, when $\gamma=0$, we have
\begin{eqnarray}
\frac{d{J}_{I}(\alpha^*,\beta^*,\gamma)}{d\gamma}\Big{|}_{\gamma=0}=-q_{I}\leq0,\nonumber
\end{eqnarray}
and thus, the maximum of ${J}_{I}(\alpha^*,\beta^*,\gamma)$ is $\gamma^*=0$ or $\gamma^*=1$. We judge the sizes of ${J}_{I}(\alpha^*,\beta^*,0)$ and ${J}_{I}(\alpha^*,\beta^*,1)$ by
\begin{eqnarray}
{J}_{I}(\alpha^*,\beta^*,1)=J_{I}(\alpha^*,\beta^*,0)+(\frac{\beta^*}{\alpha^*+\beta^*})^2-\frac{1}{2}-q_{I}.\nonumber
\end{eqnarray}
Hence, if $(\frac{\beta^*}{\alpha^*+\beta^*})^2-q_{I}<\frac{1}{2}$, the maximum of ${J}_{I}(\alpha^*,\beta^*,\gamma)$ is $\gamma^*=0$. If $(\frac{\beta^*}{\alpha^*+\beta^*})^2-q_{I}>\frac{1}{2}$, then $\gamma^*=1$. If $(\frac{\beta^*}{\alpha^*+\beta^*})^2-q_{I}=\frac{1}{2}$, then $\gamma^*=1~\text{or}~0$.

\section{Proof of Theorem 1}\label{Theorem1}
\subsection{The APT-I game with a discovered malicious insider}\label{proofA}
Based on the best response strategies in \eqref{A1}-\eqref{A3}, there are twelve possible Nash equilibria in the APT-I game. We take a look at each of them and find out feasible Nash equilibria.

\textbf{Case 1:} $\text{NE}_{1}=(\frac{r_{A}}{\beta^*},\frac{r_{D}}{\alpha^*},0).$ 

If $r_{A}\neq r_{D}$, then $\text{NE}_{1}$ does not exist, otherwise
\begin{eqnarray}
	\frac{r_{A}}{\beta^*}\leq1,~\frac{r_{D}}{\alpha^*}\leq1,~\text{and}~ (\frac{\beta^*}{\alpha^*+\beta^*})^2-q_{I}\leq\frac{1}{2}.\nonumber
\end{eqnarray}
Since $\alpha^*=\frac{r_{A}}{\beta^*}$, 
\begin{eqnarray}
	1\leq \frac{1}{\beta^*}\leq \frac{1}{r_{A}} ~\text{and}~\frac{1}{\beta^*}\geq\frac{1}{\varepsilon}.	\label{Acase11}
\end{eqnarray}
Since $q_{I}\geq 0$ and $r_{D}=r_{A}\leq 1$, $\varepsilon>r_{A}$. We obtain the following two cases:
\begin{itemize}
	\item If $q_{I}>(\frac{1}{r_{A}+1})^2-\frac{1}{2}$, then $\frac{1}{\varepsilon}<1$, which means that $\text{NE}_{1}=(\frac{r_{D}}{\beta^*},\beta^*,0)$ with $\beta^*\in[r_{D},1]$. 
	
	\item If $q_{I}\leq(\frac{1}{r_{A}+1})^2-\frac{1}{2}$, then $\frac{1}{\varepsilon}\geq1$, which means that $\text{NE}_{1}=(\frac{r_{D}}{\beta^*},\beta^*,0)$ with $\beta^*\in[r_{D},\varepsilon].$
\end{itemize}

%To sum up, We can obtain two kinds of Nash equilibrium in Case 1 as follows,
%\begin{itemize}
%\item if $r_{A}=r_{D}\leq1$ and $q_{I}>(\frac{1}{r_{A}+1})^2-\frac{1}{2}$, the Nash equilibrium falls into $\Theta_{1}=(\frac{r_{A}}{\beta^*},\beta^*,0)$ with $\beta^*\in [r_{A},1]$;
%
%\item If $r_{A}=r_{D}\leq1$ and $q_{I}\leq(\frac{1}{r_{A}+1})^2-\frac{1}{2}$, the Nash equilibrium falls into $\Theta_{1}=(\frac{r_{A}}{\beta^*},\beta^*,0)$ with $\beta^*\in [r_{A},\varepsilon]$.
%\end{itemize}

\textbf{Case 2:} $\text{NE}_{2}=(\alpha^*,\beta^*,\gamma^*)=(r_{A},1,0).$ 

If $r_{A}\geq r_{D}$, then $\text{NE}_{2}$ does not exist, otherwise
\begin{eqnarray}
	r_{A}\leq1,~r_{A}<r_{D},~\text{and}~(\frac{1}{r_{A}+1})^2-q_{I}\leq\frac{1}{2}.\label{caseA4}
\end{eqnarray}
Thus, if $r_{A}\leq1,~r_{A}<r_{D}~\text{and}~ q_{I}\geq(\frac{1}{r_{A}+1})^2-\frac{1}{2}$, then $\text{NE}_{2}=(r_{A},1,0)$.

\textbf{Case 3:} $\text{NE}_{3}=(\alpha^*,\beta^*,\gamma^*)=(1,r_{D},0).$ 

If $r_{A}\leq r_{D}$, then $\text{NE}_{3}$ does not exist, otherwise
\begin{eqnarray}
	r_{D}\leq1,~r_{D}<r_{A},~\text{and}~ (\frac{p_{D}}{p_{D}+q_{D}})^2-q_{I}\leq\frac{1}{2}.\nonumber
\end{eqnarray}
Since $p_{D}\leq q_{D}$, $q_{I}\geq(\frac{p_{D}}{p_{D}+q_{D}})^2-\frac{1}{2}$. Hence, if $r_{D}\leq1~\text{and}~r_{D}<r_{A}$, then $\text{NE}_{3}=(1,r_{D},0)$.

\textbf{Case 4:} $\text{NE}_{4}=(\alpha^*,\beta^*,\gamma^*)=(1,1,0).$ 

If $r_{A}\leq1$ or $r_{D}\leq1$, $\text{NE}_{4}$ does not exist. If $r_{A}>1$ and $r_{D}>1$, then $\text{NE}_{4}=(1,1,0)$.

\textbf{Case 5:} $\text{NE}_{5}=(\frac{r_{A}}{p_{A}\beta^*},\frac{r_{D}}{p_{D}\alpha^*},1).$ 

If $\frac{r_{A}}{p_{A}}\neq\frac{r_{D}}{p_{D}}$, then $\text{NE}_{5}$ does not exist, otherwise
\begin{eqnarray}
	\frac{r_{D}}{p_{D}\beta^*}\leq1,~\frac{r_{D}}{p_{D}\alpha^*}\leq1,~\text{and}~(\frac{\beta^*}{\alpha^*+\beta^*})^2-q_{I}\geq\frac{1}{2}.\label{caseD2}
\end{eqnarray}

Putting $\alpha^*=\frac{r_{D}}{p_{D}\beta^*}$  into \eqref{caseD2} yields
\begin{eqnarray}
	1\leq\frac{1}{\beta^*}\leq \frac{p_{D}}{r_{D}},~\frac{1}{\beta^*}\leq\frac{1}{\varepsilon\sqrt{p_{D}}},~\text{and}~q_{I}\leq (\frac{1}{\frac{r_{D}}{p_{D}(\beta^*)^2}+1})^2-\frac{1}{2}.\nonumber
\end{eqnarray}
Due to $q_{I}\geq 0$, we have $\frac{1}{\varepsilon\sqrt{p_{D}}}<\frac{p_{D}}{r_{D}}$ and the following two cases. 
\begin{itemize}
	\item If $q_{I}>(\frac{1}{\frac{r_{D}}{p_{D}}+1})^2-\frac{1}{2}$, then $\frac{1}{\varepsilon\sqrt{p_{D}}}<1$, which means that $\text{NE}_{5}$ does not exist.
	
	\item If $q_{I}\leq(\frac{1}{\frac{r_{D}}{p_{D}}+1})^2-\frac{1}{2}$, then $\frac{1}{\varepsilon\sqrt{p_{D}}}\geq1$, which means that $\text{NE}_{5}=(\frac{r_{D}}{p_{D}\beta^*},\beta^*,1)$ with $\beta^*\in[\varepsilon\sqrt{p_{D}},1]$.  
\end{itemize}

\textbf{Case 6:} $\text{\text{NE}}_{6}=(\frac{p_{A}(1-\gamma^*)^2+(\gamma^*)^2}{q_{A}\beta^*}, \frac{p_{D}(1-\gamma^*)^2+(\gamma^*)^2}{q_{D}\beta^*},\gamma^*)$ with $\gamma^*=1~\text{or}~0$. In this case, \eqref{Acase11} and \eqref{caseD2} are required to be satisfied simultaneously. Since $p_{D}\neq 1$ and $p_{A}\neq 1$, $\text{\text{NE}}_{6}$ does not exist.

\textbf{Case 7:} $\text{NE}_{7}=(\alpha^*,\beta^*,\gamma^*)=(\frac{r_{A}}{p_{A}},1,1).$ 

If $\frac{r_{A}}{p_{A}}\geq\frac{r_{D}}{p_{D}}$, then $\text{NE}_{7}$ does not exist, otherwise
\begin{eqnarray}
	\frac{r_{A}}{p_{A}}\leq1,~\frac{r_{A}}{p_{A}}<\frac{r_{D}}{p_{D}},~\text{and}~(\frac{q_{A}}{1+q_{A}})^2-q_{I}\geq\frac{1}{2}.\label{caseD5}
\end{eqnarray}
To sum up, if $\frac{r_{A}}{p_{A}}\leq1,~\frac{r_{A}}{p_{A}}<\frac{r_{D}}{p_{D}},~\text{and}~ q_{I}\leq(\frac{1}{\frac{r_{A}}{p_{A}}+1})^2-\frac{1}{2}$, then $\text{NE}_{7}=(\frac{r_{A}}{p_{A}},1,1)$.

\textbf{Case 8:} $\text{NE}_{8}=(\frac{p_{A}(1-\gamma^*)^2+(\gamma^*)^2}{q_{A}\beta^*}, 1,\gamma^*)$ with $\gamma^*=0~\text{or}~1$. In this case, \eqref{caseA4} and \eqref{caseD5} are required to be satisfied simultaneously. Since $p_{D}\neq 1$ and $p_{A}\neq 1$, $\text{\text{NE}}_{8}$ does not exist.

\textbf{Case 9:} $\text{NE}_{9}=(\alpha^*,\beta^*,\gamma^*)=(1,r_{D},1).$ 

If $r_{A}\leq r_{D}$, then $\text{NE}_{9}$ does not exist, otherwise
\begin{eqnarray}
	r_{D}\leq1,~r_{D}<r_{A}~\text{and}~ (\frac{p_{D}}{p_{D}+q_{D}})^2-q_{I}\geq\frac{1}{2}.\nonumber
\end{eqnarray}
Since $p_{D}\leq q_{D}$, $q_{I}\leq(\frac{p_{D}}{p_{D}+q_{D}})^2-\frac{1}{2}$ does not hold. Hence, $\text{NE}_{9}$ does not exist.

\textbf{Case 10:} $\text{NE}_{10}=(1,r_{D},\gamma^*)$ with $\gamma^*=0~\text{or}~1$. Since $q_{I}=(\frac{p_{D}}{p_{D}+q_{D}})^{2}-\frac{1}{2}<0$, $\text{NE}_{10}$ does not exist.

\textbf{Case 11:} $\text{NE}_{11}=(\alpha^*,\beta^*,\gamma^*)=(1,1,1).$ Since $q_{I}\leq-\frac{1}{4}$, $\text{NE}_{11}$ does not exist.

\textbf{Case 12:} $\text{NE}_{12}=(1,1,\gamma^*)$ with $\gamma^*=0~\text{or}~1$. Since $q_{I}=(\frac{1}{2})^{2}-\frac{1}{2}<0$, $\text{NE}_{12}$ does not exist.

\subsection{The APT-I game with a discovered inadvertent insider}\label{proofB}
Based on the best response strategies presented in \eqref{B1}-\eqref{B3}, there are twelve possible Nash equilibria in the APT-I game. Obviously, Cases 1-4 and 9-12 are the same as those in Appendix~\ref{proofA}. Moreover, due to $p_{D}\neq 1$, the Nash equilibria in Cases 6 and 8 do not exist. Hence, we discuss the rest cases as follows. 

\textbf{Case 5:} $\text{NE}_{5}=(\frac{r_{A}}{\beta^*},\frac{r_{D}}{p_{D}\alpha^*},1).$ 

If $r_{A} \neq\frac{r_{D}}{p_{D}}$, then $\text{NE}_{5}$ does not exist, otherwise
\begin{eqnarray}
	\frac{r_{D}}{p_{D}\beta^*}\leq1,~\frac{r_{D}}{p_{D}\alpha^*}\leq1,~\text{and}~(\frac{\beta^*}{\alpha^*+\beta^*})^2-q_{I}\geq\frac{1}{2}.\nonumber
\end{eqnarray}
The rest analysis is similar to that in Case 5 of Appendix~\ref{proofA}, and thus, is omitted. If $r_{A}=\frac{r_{D}}{p_{D}}\leq1$ and $q_{I}\leq(\frac{1}{r_{A}+1})^2-\frac{1}{2}$, then $\text{NE}_{5}=(\frac{r_{D}}{p_{D}\beta^*},\beta^*,1)$ with $\beta^*\in[\varepsilon\sqrt{p_{D}},1]$.

\textbf{Case 7:} $\text{NE}_{7}=(\alpha^*,\beta^*,\gamma^*)=(r_{A},1,1).$ 

If $r_{A}\geq\frac{r_{D}}{p_{D}}$, then $\text{NE}_{7}$ does not exist, otherwise
\begin{eqnarray}
	r_{A}\leq1,~r_{A}<\frac{r_{D}}{p_{D}},~\text{and}~(\frac{1}{r_{A}+1})^2-q_{I}\geq\frac{1}{2}.\nonumber
\end{eqnarray}
Thus, if $r_{A}\leq1, r_{A}<\frac{r_{D}}{p_{D}},~\text{and}~ q_{I}\leq(\frac{1}{r_{A}+1})^2-\frac{1}{2}$, then $\text{NE}_{7}=(r_{A},1,1)$.

Summarizing Cases 5 and 7 in the above subsections~\ref{proofA} and~\ref{proofB}, Table~\ref{table1} and Theorem 1 are obtained. Meanwhile, summarizing Cases 1-4 in Appendix~\ref{proofA}, we obtain Table~\ref{table3}.

\section{Proof of Theorem 2}\label{Theorem2}
\subsection{The APT-I game with a undiscovered malicious insider}\label{proofC}
Based on the best response strategies presented in \eqref{C1}-\eqref{C3}, there are twelve possible Nash equilibria in the APT-I game. In fact, Cases 1-4 and 9-12 are the same as those in Appendix~\ref{proofA}, and thus, are omitted. In addition, due to $p_{A}\neq 1$, the Nash equilibria in Cases 6 and 8 do not exist. We discuss the rest cases as follows. 

\textbf{Case 5:} $\text{NE}_{5}=(\frac{r_{A}}{p_{A}\beta^*},\frac{r_{D}}{\alpha^*},1).$ 

If $\frac{r_{A}}{p_{A}}\neq r_{D}$, then $\text{NE}_{5}$ does not exist, otherwise
\begin{eqnarray}
	\frac{r_{D}}{\beta^*}\leq1,~\frac{r_{D}}{\alpha^*}\leq1,~\text{and}~(\frac{\beta^*}{\alpha^*+\beta^*})^2-q_{I}\geq\frac{1}{2}.\nonumber
\end{eqnarray}
Since $\alpha^*=\frac{r_{D}}{\beta^*}$,
\begin{eqnarray}
	1 \leq \frac{1}{\beta^*} \leq \frac{1}{r_{D}}~\text{and}~\frac{1}{\beta^*}\leq\frac{1}{\varepsilon}.\nonumber
\end{eqnarray}
Since $q_{I} \geq 0$, $\varepsilon>r_{D}$. We obtain the following two cases.

\begin{itemize}
	\item If $q_{I}>(\frac{1}{r_{D}+1})^2-\frac{1}{2}$, then $\frac{1}{\varepsilon}<1$, which means that $\text{NE}_{5}$ does not exist.
	
	\item If $q_{I}\leq(\frac{1}{r_{D}+1})^2-\frac{1}{2}$, then $\frac{1}{\varepsilon}\geq1$, which means that $\text{NE}_{5}=(\frac{r_{D}}{\beta^*},\beta^*,1)$ with $\beta^*\in[\varepsilon,1]$. 
\end{itemize}

\textbf{Case 7:} $\text{NE}_{7}=(\alpha^*,\beta^*,\gamma^*)=(\frac{r_{A}}{p_{A}},1,1).$ 

If $\frac{r_{A}}{p_{A}}\geq r_{D}$, then $\text{NE}_{7}$ does not exist, otherwise
\begin{eqnarray}
	\frac{r_{A}}{p_{A}}\leq1,~\frac{r_{A}}{p_{A}}<r_{D},~\text{and}~(\frac{1}{\frac{r_{A}}{p_{A}}+1})^2-q_{I}\geq\frac{1}{2}.\nonumber
\end{eqnarray}
Thus, if $\frac{r_{A}}{p_{A}}\leq1,~\frac{r_{A}}{p_{A}}<r_{D},~\text{and}~ q_{I}\leq(\frac{1}{\frac{r_{A}}{p_{A}}+1})^2-\frac{1}{2}$, then $\text{NE}_{7}=(\frac{r_{A}}{p_{A}},1,1)$.	

\subsection{APT-I game with an undiscovered inadvertent insider}\label{proofD}
Based on the best response strategies in \eqref{D1}-\eqref{D3}, there are twelve possible Nash equilibria in the APT-I game. Similarly, Cases 1-4 and 9-12 are the same as those in Appendix~\ref{proofA}. Hence, we discuss the rest cases as follows.

\textbf{Case 5:} $\text{NE}_{5}=(\frac{r_{A}}{\beta^*},\frac{r_{D}}{\alpha^*},1).$ 

If $r_{A}\neq r_{D}$, then $\text{NE}_{5}$ does not exist, otherwise
\begin{eqnarray}
	\frac{r_{D}}{\beta^*}\leq1,~\frac{r_{D}}{\alpha^*}\leq1,~\text{and}~(\frac{\beta^*}{\alpha^*+\beta^*})^2-q_{I}\geq\frac{1}{2}.\label{Acase2}
\end{eqnarray}
The rest analysis is similar to that in Case 5 of Appendix~\ref{proofC}. Thus, if $r_{A}=r_{D}\leq1$ and $q_{I}\leq(\frac{1}{r_{D}+1})^2-\frac{1}{2}$, then $\text{NE}_{5}=(\frac{r_{D}}{\beta^*},\beta^*,1)$ with $\beta^*\in[\varepsilon,1]$.

\textbf{Case 6:} $\text{\text{NE}}_{6}=(\frac{r_{A}}{\beta^*}, \frac{r_{D}}{\alpha^*},\gamma^*)$ with $\gamma^*=1~\text{or}~0$.

In this case, \eqref{Acase11} and \eqref{Acase2} are required to be satisfied simultaneously. Thus, if $r_{A}=r_{D}\leq1~\text{and}~q_{I}=(\frac{1}{r_{D}+1})^2-\frac{1}{2}$, then $\text{NE}_{6}=(\frac{r_{D}}{\varepsilon},\varepsilon,1~\text{or}~0)$. Case 6 is a special case of Cases 1 and 5, and thus, it can be included in both of them simultaneously.

\textbf{Case 7:} $\text{NE}_{7}=(\alpha^*,\beta^*,\gamma^*)=(r_{A},1,1).$ 

If $r_{A}\geq r_{D}$, then $\text{NE}_{7}$ does not exist, otherwise
	\begin{eqnarray}
		r_{A}\leq1, r_{A}<r_{D},~\text{and}~(\frac{1}{r_{A}+1})^2-q_{I}\geq\frac{1}{2}.\label{Acase5}
	\end{eqnarray}
Thus, if $r_{A}\leq1,~r_{A}<r_{D}~\text{and}~ q_{I}\leq(\frac{1}{r_{A}+1})^2-\frac{1}{2}$, then $\text{NE}_{7}=(r_{A},1,1)$.

\textbf{Case 8:} $\text{NE}_{8}=(r_{A}, 1,\gamma^*)$ with $\gamma^*=0~\text{or}~1$.

In this case, \eqref{caseA4} and \eqref{Acase5} are required to be satisfied simultaneously. Thus, if $r_{A}\leq1,~r_{A}<r_{D}~\text{and}~q_{I}=(\frac{1}{r_{A}+1})^{2}-\frac{1}{2}$, then $\text{NE}_{8}=(r_{A}, 1,\gamma^*)$ with $\gamma^*=1~\text{or}~0$. Case 8 is a special case of Cases 2 and 7, and thus, it can be included in both of them simultaneously.

Summarizing Cases 5 and 7 in the above subsections~\ref{proofC} and~\ref{proofD}, Table~\ref{table2} and Theorem 2 can be obtained. Meanwhile, summarizing Cases 1-4 in Appendix~\ref{proofA}, we obtain Table~\ref{table3}.

\bibliographystyle{ieeetr}
\bibliography{nonconvexquantization}


%\end{IEEEbiography}


%\newpage



%\vfill


\end{document}